%% file: gw-fp-et.tex
\documentclass[11pt]{article} 

\usepackage{amsfonts}
\usepackage{latexsym}
\input{psfig-97}
\psfiginit

\font\tinybbfont=msbm6
\font\scriptsizebbfont=msbm7 scaled \magstep 1
 1
\font\footnotesizebbfont=msbm9 scaled \magstep 0
 2
\font\bbfont=msbm9 scaled \magstep1  
 1
 2
 3
 4

\def\tinyBbb#1{\hbox{\tinybbfont #1}}
\def\scriptsizeBbb#1{\hbox{\scriptsizebbfont #1}}

\def\footnotesizeBbb#1{\hbox{\footnotesizebbfont #1}}

\def\Bbb#1{\hbox{\bbfont #1}}

\input{newcmd.tex}

\begin{document}

\enlargethispage{23cm}

\begin{titlepage}

$ $

\vspace{-2cm} 

\noindent\hspace{-1cm}
\parbox{6cm}{\small April 2005}\
   \hspace{6.5cm}\
   \parbox{5cm}{math.AG/0505084}

\vspace{1.5cm}

\centerline{\large\bf
 Transformation of algebraic Gromov-Witten invariants
  of three-folds
}
\vspace{1ex}
\centerline{\large\bf
 under flops and small extremal transitions,
}
\vspace{1ex}
\centerline{\large\bf
 with an appendix from the stringy and the symplectic viewpoint
}

\vspace{1.5cm}
\centerline{\large
  Chien-Hao Liu
  \hspace{1ex} and \hspace{1ex}
  Shing-Tung Yau
}

\vspace{3em}

\begin{quotation}
\centerline{\bf Abstract}
\vspace{0.3cm}
\baselineskip 12pt  
{\small
 We study how Gromov-Witten invariants of projective $3$-folds
  transform under a standard flop and a small extremal transition
  in the algebro-geometric setting from the recent development
  of algebraic relative Gromov-Witten theory and its applications.
 This gives an algebro-geometric account
  of Witten's wall-crossing formula for correlation functions
   of the descendant nonlinear sigma model in adjacent geometric
   phases of a gauged linear sigma model and
  of the symplectic approach in an earlier work of An-Min Li and
   Yongbin Ruan on the same problem.
 A terse account from the stringy and the symplectic viewpoint
  is given in the appendix to complement and compare to
  the discussion in the main text.
} 
\end{quotation}

\vspace{12em}

\baselineskip 12pt
{\footnotesize
} 

\medskip

\noindent {\small
MSC number 2000$\,$:
 14N35; 53D45, 81T30.
} 

\bigskip

\baselineskip 11pt
{\footnotesize
\noindent{\bf Acknowledgements.}
 We thank
  Jun Li and Kefeng Liu
 for discussions on the project and 
  Yongbin Ruan
 for communications and a very detailed explanation on a technical
 point in [L-R], which resolves a puzzle we originally have in
 comparing the symplectic and the algebraic statements for
 GW-invariants under small extremal transitions.
 C.-H.L.\ would like to thank in addition
  Joe Harris and Mihnea Popa
 for valuable lectures,
  Shiraz Minwalla
 for the insightful lectures on QFT and its phases in spring 2004,
  Uwe Storch
 for a discussion, and
  Ling-Miao Chou for the moral support.
 The project is supported by NSF grants DMS-9803347 and DMS-0074329.
} 

\end{titlepage}

\newpage
\begin{titlepage}

$ $
\vspace{10em}

\centerline{\small\it
  C.-H.L.\ dedicates this review
 }
\centerline{\small\it
  to the numerous teachers who educated him
 }
\centerline{\small\it
 and to Ling-Miao Chou for her tremendous love.
}
\vspace{10em}

A reflection on $\,${\it string/M/F/$\cdots?\,$-$\,$theory}$\,$,
 cf.\ Sec.~A.1$\;$:

\bigskip

\begin{figure}[htbp]
 \centerline{\psfig{figure=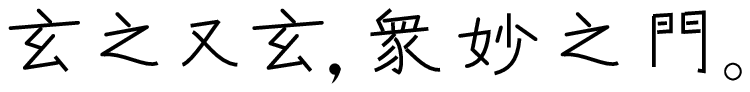,width=30pc,caption=}}
\end{figure}

\vspace{-4em}
\begin{itemize}
 \item[]
 ($\,${\it Mystery and beyond mystery, door to all magics.}$\,$)
\end{itemize}

\bigskip
\begin{figure}[htbp]
 \centerline{\psfig{figure=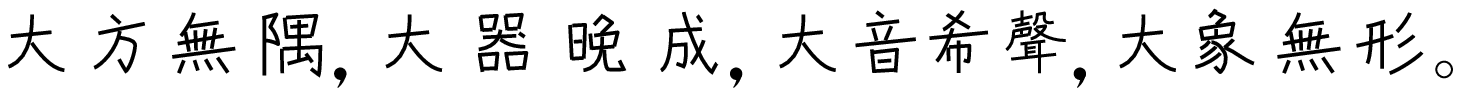,width=30pc,caption=}}
\end{figure}

\vspace{-4em}
\begin{itemize}
 \item[]
 ($\,${\it So large that it has no bounds; \newline
 so big that it takes a long time to make; \newline
 so harmonious that it fits no tunes; \newline
 so beautiful that it assumes no shapes.}$\,$)
\end{itemize}

\bigskip
\bigskip
\noindent
 \hspace{15em}$\sim$
 \parbox[t]{20em}{\small
  Lao-Tzu (600 B.C.), {\sl Tao-te Ching} \newline
  ({\sl The Scripture on the Way and its Virtue}), \newline
  excerpt from Chapter~1 and Chapter~41.  \newline\newline
  {\footnotesize English translation by Ling-Miao Chou.}
 } 

\end{titlepage}

\newpage
$ $

\vspace{-4em}  

\centerline{\sc
 Gromov-Witten Invariants under Flops and Extremal Transitions
}

\vspace{2em}

\baselineskip 14pt  

\begin{flushleft}
{\Large\bf 0. Introduction and outline.}
\end{flushleft}
Surgeries on and transitions of manifolds/varieties are important
 issues in classifications of K\"{a}hler manifolds and varieties.
Once the Gromov-Witten theory for smooth manifolds/varieties are
 understood
 (see [Be], [Si1] for reviews and
   [L-T2], [Si2] for comparisons of the symplectic and
    the algebro-geometric construction),
it is very natural that one wants to extend this understanding to
 transformations of Gromov-Witten invariants when the manifold/variety
 changes.
Indeed, this theme lies on the cross-road of
 superstring theory, symplectic geometry, and algebraic geometry.

The generating function of Gromov-Witten invariants of
 a Calabi-Yau manifold $X$ corresponds to the correlation functions
 of the A-model for $X\,$, [Wit4].
The transition of the Calabi-Yau manifolds corresponds to moving
 from one phase to another of the $d=2$, $N=(2,2)$ superconformal
 field theory (SCFT) with a specified number of chiral multiplets
 depending on the dimension of $X$;
and the associated transformation of Gromov-Witten invariants
 corresponds to the transformation of the correlation functions
 of the associated A-models in the different phases of SCFT.
In the case of flops of Calabi-Yau $3$-folds realizable
 as complete-intersection subvarieties of a toric variety,
Witten [Wit5] (year 1993) has put this picture manifestly
 as a transition between adjacent geometric phases of
 a ($d=2$, $N=(2,2)$) gauged linear sigma model
 (GLSM; [Wit5], see also [M-P], [H-V] and [H-I-V]) and
conjectured a wall-crossing formula
 (cf.\ [Wit5: Sec.~5.5 and Eq.(5.48)])
 by isolating the effect of the rigid curve involved in the flop
 to the full correlation function (expressed as an instanton sum).
The reason that this remained only a conjecture as stated in [Wit5]
 (even assuming the multiple cover formula at that time)
 is that it is not clear that the effect of the original rigid curve
 to the instanton sum of the flopped Calabi-Yau $3$-fold can be
 really isolated out in such a simple way since curves can break
 under deformations and there are also possible contributions to
 Gromov-Witten invariants from reducible curves with some but
 not all component wrapping the rigid curve.

A such technicality needed to understand the exact details was
 later fulfilled in the work of A.-M.~Li and Y.~Ruan in [L-R]
 (year 1998) within the program of a general study of Gromov-Witten
 invariants in birational geometry.
There they developed a symplectic relative Gromov-Witten theory
 and derived a gluing formula for symplectic cuts,
 (see also [I-P1] and [I-P2] of E.-N.\ Ionel and T.H.\ Parker
  for related works).
This enabled them
 to isolate the common summand and the different summand
  in the instanton sum of $3$-folds on the two sides of
  a standard flop or a small extremal transition and
 to understand their changes.
How the curves that are involved in the surgery affect
 the transformation of Gromov-Witten invariants under flops
 can be extracted by using the gluing formula as a medium;
a complete statement of Witten's formula can then be spelled out and
 justified ([L-R: Theorem A, Corollary A.2, Theorem B, Corollary B.2]).

In the history of the development of Gromov-Witten theory,
 the symplectic construction and the algebro-geometric construction
 have been going side-by-side with each other.
It is thus desirable to also have a pure algebro-geometric
 construction parallel to the symplectic relative Gromov-Witten
  theory in [L-R], [I-P1], and [I-P2].
This was realized by J.~Li in [Li1] and [Li2], who constructed
 an algebraic relative Gromov-Witten theory and derived
 a degeneration formula with respect to a relative ample line bundle
 for a projective family over a smooth curve with the degenerate
 fiber $Y_1\cup_D Y_2$ a gluing of two smooth variety-divisor pairs
 $(Y_i,D)$.
His work was later used in [L-Y1] and [L-Y2] to understand
 Gromov-Witten invariants associated to a curve class
  for smooth projective varieties under blow-ups along a smooth locus
   and
  for a projective conifold that fits in a degeneration with smooth
   total space.
These are the algebro-geometric counterparts of the main tools
 of [L-R] in symplectic category.
The remaining discussions needed to lead to the same statements
 of [L-R] is given in this note.
Subject to the fact that choices of cycles in the algebraic geometry
 are more restrictive than in symplectic geometry, this puts Witten's
 wall-crossing formula and the general results of [L-R]
 in the algebro-geometric category.

Finally, it should be noted that there are more general flip-flops
 and extremal transitions in the study of birational geometry and
 classifications of K\"{a}hler manifolds.
We conclude this introduction with a remark from Y.~Ruan
 on a comparison of the symplectic and the algebro-geometric setting
 for transformations of Gromov-Witten invariants under flops and
 small extremal transitions:

\bigskip

\noindent
{\it Remark $[$symplectic vs.\ algebraic$\,]$.} ([Ru3].)
 In algebro-geometric category there are many flops other than
  the standard flop. How to decompose them into
  a composition of blow-ups and blow-downs can be complicated.
 In contrast, in the symplectic category one can perform a local  
  holomorphic deformation to deform extremal rational curves into
  a disjoint union of rational curves of
  ${\cal O}_{{\scriptsizeBbb P}^1}(-1)
    \oplus {\cal O}_{{\scriptsizeBbb P}^1}(-1)$-type
  curve neighborhood.
 These curve neighborhood can be patched with the complex structure
  on the complement of these curves to form a tamed almost complex
  structure.
 In such a fashion, one can decompose a general flop into a disjoint
  union of standard flops in the almost complex category.
 In this way, the results of [L-R] work for any flops instead of
  just the standard flop of $3$-folds.
 Similarly for small extremal transitions.

\bigskip

\noindent
Readers are referred to [L-R: Sec.~1], [Mo2], and [Ru2]
 for an overview of Gromov-Witten theory in this direction.
It can be challenging to give a complete algebro-geometric treatment
 for Gromov-Witten invariants under general
 flips/flops/extremal-transitions as well.

\bigskip

For the contents of the note,
in Sec.~1 and Sec.~2, we recall the basic definitions and results
 from [Li1], [Li2], [L-Y1], and [L-Y2] needed for the discussions.
In Sec.~3, we derive from these the transformation formula of
 Gromov-Witten invariants of projective $3$-folds under
 a standard flop and a small extremal transition.
A terse account of the problem from the stringy and the symplectic
 viewpoint is given in the Appendix to compare to and to complement
 the algebro-geometric viewpoint in the main text.

\bigskip

\noindent
{\it Convention.}
 Standard notations, terminology, and operations in algebraic
  geometry can be found in [Ha], [Kol-M], and [Fu].
 All schemes are over ${\Bbb C}$.
 The word ``relative" has two different but complementary meanings:
  one from topology, e.g.\ 
   a {\it relative map} to a smooth variety-divisor pair $(X,D)$
    means a map to $X$ with image intersecting $D$ properly;
   and
  the other from algebraic geometry in Grothendieck's formulation
   meaning relative to a base scheme, e.g.\
    a {\it relative cycle} on $W/{\Bbb A}^1$ in this note means
     a formal linear combination of subvarieties on $W$ that
     are flat over ${\Bbb A}^1$.

\bigskip

\begin{flushleft}
{\bf Outline.}
\end{flushleft}
{\small
\baselineskip 11pt  
\begin{quote}
 1. Algebraic relative Gromov-Witten theory and degeneration
    formulas.

 2. Degenerations associated to blow-ups and conifolds.

 3. Transformation of GW-invariants of $3$-folds.
    \vspace{-1ex}
    \begin{quote}
     \hspace{-1.3em}
     3.1 \ \parbox[t]{12cm}{
     Gromov-Witten invariants of $3$-folds under a standard flop.}

     \hspace{-1.3em}
     3.2 \ \parbox[t]{12cm}{
     Gromov-Witten invariants of $3$-folds under a small extremal
     transition.}
    \end{quote}

 \vspace{-.8ex}
 Appendix. The stringy and the symplectic aspect.
    \vspace{-1ex}
    \begin{quote}
     \hspace{-1.3em}
     A.1 \ \parbox[t]{12cm}{
     Transformation of GW-invariants from the stringy viewpoint.}

     \hspace{-1.3em}
     A.2 \ \parbox[t]{12cm}{
     Transformation of GW-invariants from the symplectic viewpoint.}
    \end{quote}
\end{quote}
} 

\bigskip

\baselineskip 14pt  

\section{Algebraic relative Gromov-Witten theory and
         degeneration formulas.}

We recall the basic ingredients of the algebraic relative
 Gromov-Witten theory needed for the discussion.
See [Li1], [Li2], and [G-V] for more details.

\bigskip

\begin{flushleft}
{\bf Algebraic relative Gromov-Witten theory.}
\end{flushleft}
Let $(Y,D)$ be a smooth variety-divisor pair with $Y$ projective.
A {\it relative map} $f:C\rightarrow (Y,D)$ is a morphism from
 a prestable curve $C$ to $Y$ such that $f$ meets $D$ properly
 (i.e.\ $f^{-1}(D)$ is a divisor on $C$).
Fix a topological type $\Gamma=(g, n; \beta)$
 where $g$ is the genus of a prestable curve, $n$ is the number
  of marked points on that curve and $\beta\in H_2(Y;{\Bbb Z})$,
the moduli stack ${\cal M}_{g,n}(Y,D\,;\beta)$ of relative maps of
 topological type $\Gamma$ and with finite automorphism group 
 is a Deligne-Mumford stack, from the Hilbert-scheme construction
 in [F-P].
Its compactification in the usual moduli stack of
 $\overline{\cal M}_{g,n}(Y;\beta)$ of stable maps to $Y$ of
  topological type $\Gamma$ can contain maps $f:C\rightarrow Y$
 that send a subcurve of $C$ to $D$.
When such degeneracy occurs in an ${\Bbb A}^1$-family of
 stable maps to $Y$ with generic fiber a relative map, 
it can be removed
 by a finite base change,
     of degree bounded above by a number determined by $\Gamma$,
     and then
    refilling both the domain curve and the target space of
     the bad fibers.
This procedure extends the target space of relative maps
 in question to the fibers of the universal family of the
 {\it stack ${\mathfrak Y}^{\rel}$ of expanded relative pairs
      associated to} $(Y,D)$.
${\mathfrak Y}^{\rel}$ with its universal family is the descent
 of the local standard models of expanded relative pairs
 $(Y[m], D[m])\rightarrow {\Bbb A}^m$ constructed in
 [Li1: Sec.~4.1], (see also [G-V: Sec.~2.6]).
As an abstract stack, ${\mathfrak Y}^{\rel}$ is isomorohic to
 the open substack ${\cal T}$ of $\overline{\cal M}_{0,3}$
  that parameterizes prestable curves $({\Bbb P}^1_{[m]}; 0,1,\infty)$
  of arithmetic genus $0$ with all nodes separating $\infty$ from
  $\{0,1\}$.
The fiber of the universal family of ${\mathfrak Y}^{\rel}$
 over the ${\Bbb C}$-point $({\Bbb P}^1_{[m]}; 0,1,\infty)$ of
 ${\mathfrak Y}^{\rel}$ is given by the gluing of $Y$ with
 an $m$-chain of ruled varieties
 $\Delta_i:= {\Bbb P}({\cal O}_D\oplus {\cal N}_{Z/Y})$ over $D\,$,
 $i=1,\,\ldots,\,m\,$:
 $$
  Y_{[m]}\;
   :=\; Y\, \cup_{D=D_{1,0}} \Delta_1 \cup_{D_{1,\infty}=D_{2,0}}\;
            \cdots\; \cup_{D_{m-1,\infty}=D_{m,0}} \Delta_m
 $$
 with
 $$
   D_{[m]}\; =\; D_{m,\infty}\; \subset\; \Delta_m\,.
 $$
 Here,
  ${\cal N}_{D/Y}$ is the normal bundle of $Z$ in $Y$,
  $D_{i,0}$ (resp.\ $D_{i,\infty}$) is the section of
   $\Delta_i\rightarrow D$ corresponding to ${\cal N}_{D/Y}$
   (resp.\ ${\cal O}_D$), and
  the gluing $D_{i,\infty}=D_{i+1,0}$ is given by the projection
   $\Delta_i\rightarrow D$\,,
 cf.\ {\sc Figure} 1-1.
 %
 %
By construction, there is a tautological morphism
 $\varphi: (Y[m],D[m])\rightarrow (Y,D)$ that pinches all
 the $\Delta_i$ in fibers of $Y[m]$ over ${\Bbb A}^m$ to $D$.

\begin{figure}[htbp]
 \setcaption{\small
 {\sc Figure} 1-1.
  \baselineskip 12pt
  A local model $(Y[m], D[m])/{\Bbb A}^m$ for the universal family
   of the stack ${\mathfrak Y}^{\rel}$ of expanded relative pairs
   associated to $(Y,D)$ and its tautological morphism to $(Y,D)$.
  The fiber over the origin $0\in {\Bbb A}^m$ is $(Y_{[m]},D_{[m]})$.
  Illustrated here is the case $m=2$.
 } 
\centerline{\psfig{figure=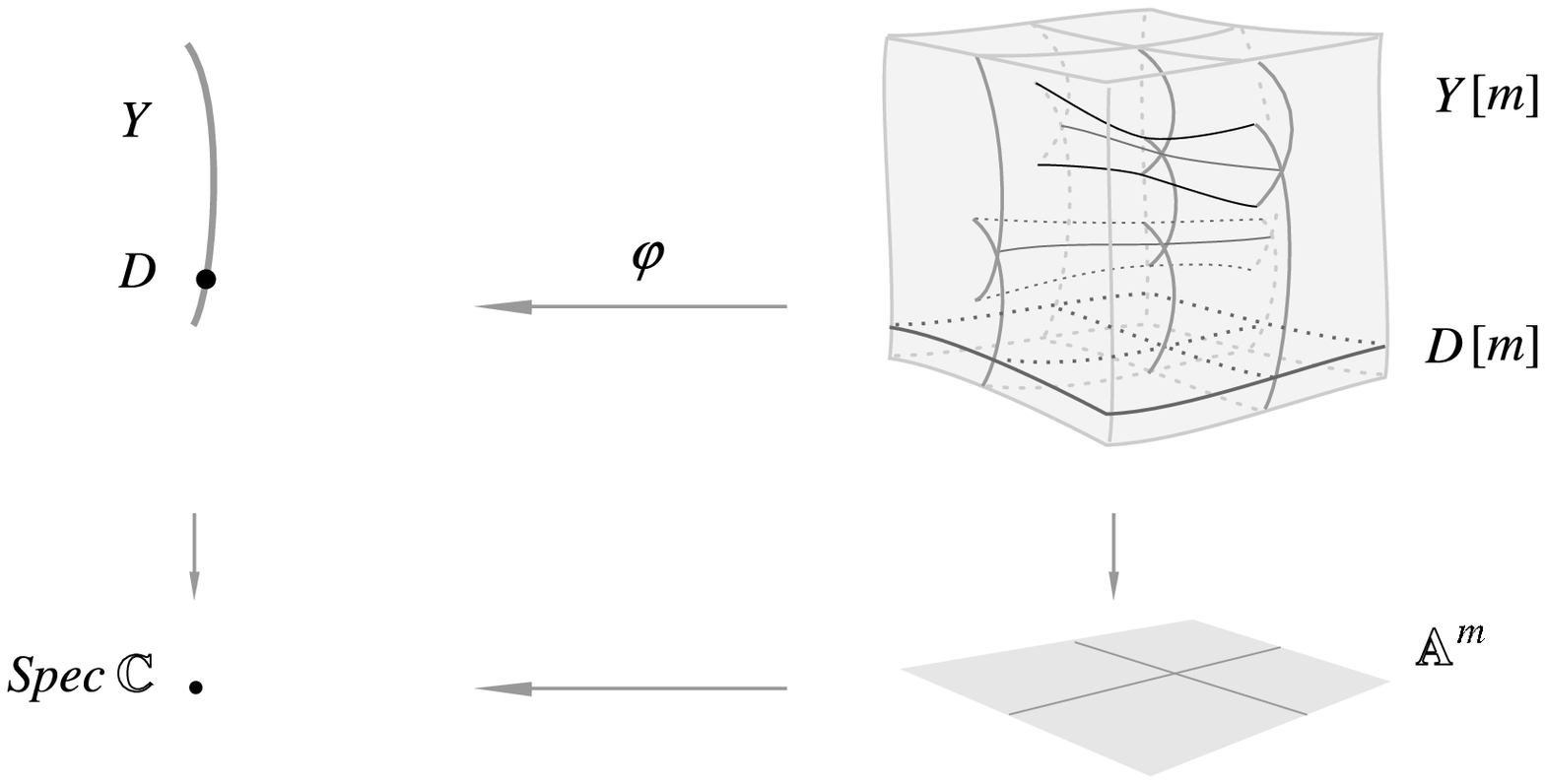,width=13cm,caption=}}
\end{figure}

\bigskip

\noindent
{\bf Definition 1.1 [admissible weighted graph].}
([Li1: Definition 4.6].)
 {\rm
  Given a relative pair $(Y,D)$, an {\it admissible weighted graph}
   $\Gamma$ for $(Y,D)$ is a graph without edges together
   with the following data:
  \begin{itemize}
   \item [(1)]
    an ordered collection of {\it legs},
    an ordered collection of weighted {\it roots}, and
    two {\it weight functions} on the vertex set
     $g:V(\Gamma)\rightarrow {\Bbb Z}_{\ge 0}$ and
    $b:V(\Gamma)\rightarrow H_2(Y;{\Bbb Z})$,

   \item[(2)]
    $\Gamma$ is {\it relatively connected} in the sense that either
     $|V(\Gamma)|=1$ or each vertex in $V(\Gamma)$ has at least one root
     attached to it.
  \end{itemize}
} 

\bigskip

Given an admissible weighted graph $\Gamma$ for $(Y,D)$
 with $l$ vertices, $n$ legs, and $r$ roots,
a {\it stable relative map} of type $\Gamma$ to the fibers of
 the universal family of ${\mathfrak Y}^{\rel}$ is locally
 modelled on a diagram:
 (in notation, $f:{\cal C}/S\rightarrow Y[m]/{\Bbb A}^m$)
 $$
  \begin{array}{ccc}
   {\cal C}    & \stackrel{f}{\longrightarrow}   & Y[m]            \\
   \downarrow  &                                 & \downarrow      \\
   S           &  \rightarrow                    & {\Bbb A}^m
  \end{array}\,,
 $$
 where
 \begin{itemize}
  \item[$\cdot$]
   ${\cal C}$ is a disjoint union of flat families ${\cal C}_i$,
    $i=1,\,\ldots,\,l\,$, of prestable curves over $S$ of arithmetic
    genus $g(v_i)$;

  \item[$\cdot$]
   $p_i:S\rightarrow {\cal C}$, $i=1,\,\ldots,\,n$, and
   $q_j:S\rightarrow {\cal C}$, $j=1,\,\ldots,\,r$ are disjoint
    sections away from the singular locus of fibers of ${\cal C}$
    over $S$
   such that $p_i(S)\subset {\cal C}_k$
     (resp.\ $q_j(S)\subset {\cal C}_k$)
     if $i$-th leg (resp.\ $j$-th root) is attached to
     the vertex $v_k$ of $\Gamma$;

  \item[$\cdot$]
   $f$ is a nondegenerate predeformable morphism to $Y[m]$
    such that
     $f^{-1}(D[m])=\sum_{j=1}^r\mu_j q_j(S)$,
       where $\mu_j$ is the weight associated to the $j$-th root,
    and that
     at each closed point $s\in S$,
     $\varphi_{\ast}(f({{\cal C}_i}_s))=b(v_i)$
     in $H_2(Y;{\Bbb Z})$;   

  \item[$\cdot$]
   at each closed point $s\in S$, $f_s$ has only a finite
    automorphism group,
 \end{itemize}
([Li1: Sec.~2.2, Sec.~3.1, Sec.~4.2]).
We remark that in the above definition, an automorphsim of
 $f:C\rightarrow Y_{[m]}$ is a pair
 $(a,b)\in \Aut(C)\times Aut(Y_{[m]}/Y)\,$
 (i.e.\ $b$ is an automorphism on $Y_{[m]}$ that descends
        to the identity map on $Y$ via $\varphi\,$)
 such that $f\circ a=b\circ f\,$,
(cf.\ [Li1: beginning of Sec.\ 3.1] for the definition of
 isomorphism classes and automorphisms of relative maps
 to an extended relative pair, see also [G-V: Sec.\ 2]).
The space of all such stable relative maps of type $\Gamma$
 is a Deligne-Mumford stack
 ${\mathfrak M}(Y,D\,;\Gamma)$
 ($={\mathfrak M}({\mathfrak Y}^{\rel},\Gamma)$ in
  [Li1], [Li2], [L-Y1], [L-Y2];
  here we adopt a notation change to make what is involved explicit),
[Li1: Definition 4.9, Theorem 4.10].

Suppose that $\Gamma$ has $n$ legs and $r$ roots,
 then there is an evaluation map
 $$
  \ev\,:\; {\mathfrak M}(Y,D\,; \Gamma)\;
     \longrightarrow\; Y^n
 $$
 associated to the ordinary $n$ marked points on the domain curve
 of a relative stable morphism and a distinguished evaluation map
 $$
  {\mathbf q}\;:\; {\mathfrak M}(Y,D\,; \Gamma)\;
     \longrightarrow\; D^{\,r}
 $$
 associated to the $r$ distinguished marked points that are required
 to be the only points that are mapped to $D[m]$ in a local model.

The obstruction theory associated to the deformation problems related
 to the stack ${\mathfrak M}(Y,D\,;\Gamma)$ are studied
 in [Li2: Sec.~1 and Sec.~5], see also [G-V: Sec.~2].
A cohomological description of the deformations of the separate
 constituents of a stable relative map and the natural clutching
 morphisms that relate the various separate deformation of
 the constituents are given there.
The perfectness of the obstruction theory and hence virtual
 fundamental class on
  $[{\mathfrak M}(Y,D\,; \Gamma)]^{\virt}$
 are proved and constructed in [Li2: Sec.~2].
(See also [G-V: Sec.~2.8 and Sec.~2.9].)
The relative Gromov-Witten invariants of the pair $(Y,D)$
 associated to $\Gamma$ are then defined by
 $$
  \begin{array}{cccccl}
   \Psi^{(Y,D)}_{\Gamma} & :
    & H^{\ast}(Y)^{\times n} \times H^{\ast}({\mathfrak M}_{\Gamma})
    & \longrightarrow  & H_{\ast}(D^r) \\[.6ex]
   && (\,\alpha\;,\;\zeta\,)
    & \longmapsto
    & {\mathbf q}_{\ast}
         \left(\, \rule{0em}{1.2em}
           \ev^{\ast}(\alpha)
            \cup \pi_{\Gamma}^{\ast}(\zeta)\,
                  [{\mathfrak M}(Y,D\,; \Gamma)]^{\virt}\,
         \right)   &,
  \end{array}
 $$
 where ${\mathfrak M}_{\Gamma}$ is the moduli stack of
  (possibly disconnected) stable nodal curves of the topological
  type specified by $\Gamma$ and
 $\pi_{\Gamma}:
  {\mathfrak M}(Y,D\,;\Gamma)\rightarrow{\mathfrak M}_{\Gamma}$
 is the forgetful morphism.

\bigskip

\begin{flushleft}
{\bf The degeneration formula with respect to a relative ample
     line bundle.}
\end{flushleft}
Let $\pi:(W,W_0)\rightarrow ({\Bbb A}^1,{\mathbf 0})$
 be a degeneration with
  the total space $W$ smooth,
  $\pi$ projective,
  the fiber $W_t$ over $t\ne {\mathbf 0}$ smooth, and
  the fiber $W_0$ over ${\mathbf 0}\in {\Bbb A}^1$
   the gluing  $Y_1\cup_D Y_2$ of smooth varieties $Y_i$
    along isomorphic smooth divisor $D_i\simeq D$.
Fix a relative ample line bundle $H$ on $W/{\Bbb A}^1$.
Let $(g,n;d)$ be a triple of integers.
Then a moduli stack ${\mathfrak M}({\mathfrak W},(g,n;d))$
 of stable maps from (connected) prestable curves of topological
 type $(g,n)$
 ($\,= (\mbox{arithmetic}, \mbox{number of marked points})\,$)
 to fibers of $\pi$ with $H$-degree $d$
 is constructed in [Li1].
By construction, ${\mathfrak M}({\mathfrak W}, (g,n;d))$
 fibers over ${\Bbb A}^1$ such that the fiber
 ${\mathfrak M}({\mathfrak W}, (g,n;d))_t$ over $t\ne {\mathbf 0}$
 is (a disjoint union over curve classes of $H$-degree $d$ of)
 the usual moduli stack of stable maps to the smooth $W_t$
while the fiber ${\mathfrak M}({\mathfrak W}, (g,n;d))_{\mathbf 0}$
 over ${\mathbf 0}$ gives the moduli stack of stable maps to
 the singular fiber $W_0$,
 (which is new in algebro-geometric category).

The deformation-obstruction theory of this moduli problem
 and the perfectness of the tangent-obstruction complex on
 ${\mathfrak M}({\mathfrak W},(g,n;d))$ is studied in [Li2].
The virtual fundamental class
 $[{\mathfrak M}({\mathfrak W},(g,n;d))]^{\virt}$ thus constructed
 fibers over ${\Bbb A}^1$ as well, with the property that
 its restriction to each fiber ${\mathfrak M}({\mathfrak W},(g,n;d))_t$
 gives the virtual fundamental class 
 $[{\mathfrak M}({\mathfrak W},(g,n;d))_t]^{\virt}$
 identical to the one constructed directly on
 the stack ${\mathfrak M}({\mathfrak W},(g,n;d))_t$.
(With appropriate cycles inserted via evaluation maps),
 this gives a constancy over ${\Bbb A}^1$ of Gromov-Witten invariants
 of fibers $W_t$ of $W\rightarrow {\Bbb A}^1$ for all $t$.
The Gromov-Witten invariants of $W_0=Y_1\cup_D Y_2$ can be 
 further expressed in terms of the relative Gromov-Witten invariants
 of $(Y_1,D)$ and $(Y_2,D)$.
The upshot is a degeneration formula that expresses
 a summation of the usual Gromov-Witten invariants of a fiber
 $W_t$, $t\ne {\mathbf 0}$ in terms of a combination
 of relative Gromov-Witten invariants of pairs $(Y_1,D)$ and
 $(Y_2,D)\,$.

To state the degeneration formula precisely, recall the following
 definition:

\bigskip

\noindent
{\bf Definition 1.2 [admissible triple].}
([Li1: Definition 4.11].)
{\rm
 Given a gluing $Y_1\cup_D Y_2$ of relative pairs,
 let $\Gamma_1$ and $\Gamma_2$ be a pair of admissible weighted graphs
  for $(Y_1,D)$ and $(Y_2,D)$ respectively.
 Suppose that $\Gamma_1$ and $\Gamma_2$ have identical number $r$
  of roots and $n_1$-many and $n_2$-many legs respectively.
 Let $n=n_1+n_2$ and $I\subset \{1,\,\ldots,\, n\}$ be a set of $n_1$
  elements.
 Then $(\Gamma_1,\Gamma_2,I)$ is called an {\it admissible triple} if
  the following conditions hold:
  \begin{itemize}
   \item[(1)]
    the weights on the roots of $\Gamma_1$ and $\Gamma_2$ coincide:
     $\mu_{1,i}=\mu_{2,i}$, $i=1,\,\ldots,\, r\,$;

   \item[(2)]
    after connecting the $i$-th root of $\Gamma_1$ and the $i$-th root
    of $\Gamma_2$ for all $i$, the resulting new graph with $n$ legs
    and no roots is connected.
  \end{itemize}
} 

\bigskip

\noindent
Given an admissible triple $\eta=(\Gamma_1,\Gamma_2,I)$ as above
 with $Y_1\cup _D Y_2=$ the degenerate fiber $W_0$ of $W/{\Bbb A}^1$,
 one has the genus function
 $$
  g(\eta)\;
  :=\;  r+1-|V(\Gamma_1\,
               \raisebox{.2ex}{\scriptsize $\coprod$}\, \Gamma_2)|
         + \sum_{v\in V(\Gamma_1)\cup V(\Gamma_2)}\,g(v)
 $$
 and the $H$-degree function
 $$
  d(\eta)\;
   :=\;
     \sum_{v\in V(\Gamma_1)}\,b_{\Gamma_1}(v)\cdot H|_{Y_1}
      + \sum_{v\in V(\Gamma_2)}\, b_{\Gamma_2}(v)\cdot H|_{Y_2}\,.
 $$
Denote by $|\eta|$ the triple of integers
 $(g(\eta), n=n_1+n_2; d(\eta))$.
For each $\eta$, one has a gluing morphism
 $$
  \Phi_{\eta}:
   {\mathfrak M}({\mathfrak Y}_1^{\it rel},\Gamma_1)
    \times_{E^r} {\mathfrak M}({\mathfrak Y}_2^{\rel},\Gamma_2)
  \rightarrow {\mathfrak M}({\mathfrak W},(g,n;d))\,,
 $$
 which is finite \'{e}tale of pure degree $|\Eq(\eta)|$ to
  its image
   ${\mathfrak M}
    ({\mathfrak Y}_1^{\rel}\,\sqcup\,{\mathfrak Y}_2^{\rel},\eta)$
  in and topologically isomorphic to
  ${\mathfrak M}({\mathfrak W}_0,\eta)$.
 [Li1: Sec.\ 4.2] and [Li2: Sec.\ 3.2].
(Here $\Eq(\eta)$ is the set of permutations of the $r$-many roots in
 $\Gamma_1$ that leaves $\eta$ unchanged.)

Given $(g,n;d)$, let
 $\Omega_{(g,n;d)}^H$ be the set of admissible triples $\eta$
  for the gluing $Y_1\cup _D Y_2$ such that $|\eta|=(g,n;d)$,
 $\overline{\Omega}_{(g,n;d)}^H$ be the set of equivalence classes
  in $\Omega_{(g,n;d)}^H$ from re-ordering of roots, and
 ${\mathbf m}(\eta)$ be the product of the weight of roots of
  $\Gamma_1$ in $\eta\in\overline{\Omega}_{(g,n;d)}^H$.
 For $\eta\in\overline{\Omega}_{(g,n;d)}^H$, assume that
  $G_{\eta}^{\ast}(\zeta)
    =\sum_{j\in K_{\eta}}\zeta_{\eta,1,j}\bboxtimes \zeta_{\eta,2,j}$,
  where
  $G_{\eta}:
         {\mathfrak M}_{\Gamma_1}\times {\mathfrak M}_{\Gamma_2}
              \rightarrow {\mathfrak M}_{g,n}$
  is the natural morphism between the related moduli stacks
  of nodal curves.
Then J.\ Li's degeneration formula of Gromov-Witten invariants
 in numerical form reads: ([Li2: Sec.\ 3 and Sec.\ 4])
 \begin{eqnarray*}
  \lefteqn{
   \Psi_{(g,n;d)}^{W_t}(\alpha(t),\zeta) } \\[.6ex]
   &&
    =\; \sum_{\eta\in\overline{\Omega}_{(g,n;d)}^H}\;
      \frac{{\mathbf m}(\eta)}{|\Eq(\eta)|}\,
        \sum_{j\in K_{\eta}}\,
         \left[ \Psi_{\Gamma_1}^{(Y_1,D)}(j_1^{\ast}\alpha(0),
                  \zeta_{\eta,1,j})\,
                \bullet\,
                \Psi_{\Gamma_2}^{(Y_2,D)}(j_2^{\ast}\alpha(0),
                  \zeta_{\eta,2,j})
         \right]_0\,,
 \end{eqnarray*}
 where
  $\Psi_{(g,n;d)}^{W_t}(\alpha(t),\zeta)$, $t\ne{\mathbf 0}$,
   is the summation of the usual algebraic Gromov-Witten invariants
   of $W_t$ over curves classes on $W_t$ of $H$-degree $d\,$;
  $\alpha\in H^0_c(R^{\ast}\pi_{3\ast}{\Bbb Q}_{W})^{\times n}\,$;
  $j_i$ is the inclusion map $Y_i\hookrightarrow W_0\,$; 
  $\bullet$ the intersection product in $H_{\ast}(D^r)$,
   where $r$ is the number of roots in either of the admissible
   graphs in the admissible triple $\eta\,$; and
  $[\;\cdot\;]_0$ is the degree-$0$ component of elements
   in $H_{\ast}(D^r)$.

\bigskip

\section{Degenerations associated to blow-ups and conifolds.}

For our applications, the degeneration formula in Sec.~1
 has to be refined to one with respect to curve classes.
This was done in [L-Y1] and [L-Y2] for the two situations:
 blow-up along a smooth subvariety and
 the conifold degeneration of $3$-folds.
These are closely related to flops and small extremal transitions
 of $3$-fold.
We recall the main results to be used in Sec.~3 and leave readers
 to [L-Y1] and [L-Y2] for more details.

\bigskip

\noindent
{\it Remark.\ 2.1 $[$homological vs.\ numerical equivalence$]$.}
 For a smooth projective variety $X$ over ${\Bbb C}$,
  the group of homological equivalence classes of (algebraic)
  curve classes in $H_2(X;{\Bbb Z})$ embeds in the group $N_1(X)$
  of numerical equivalence classes of (complex) $1$-cycles.
 Since our focus is on curve classes, the discussion and statements
  in [L-Y1] and [L-Y2] via $N_1(X)$ can be converted to ones on
  the set of curve classes in $H_2(X;{\Bbb Z})$. ([Fu] and [Lie].)
 For the same reason, all the torsion classes in $H_2(X;{\Bbb Z})$
  can be ignored.

\bigskip

\begin{flushleft}
{\bf Blow-ups.}
\end{flushleft}
Let $X$ be a smooth projective variety and $Z\subset X$
 be a smooth subvariety of $X$.
Associated to $(X,Z)$ is the degeneration
 $\pi:W\rightarrow X\times {\Bbb A}^1$ from the blow-up of
 $X\times{\Bbb A}^1$ along $Z\times{\mathbf 0}$,
 where ${\mathbf 0}$ is a point on ${\Bbb A}^1$.
By construction, $W/{\Bbb A}^1$ has the degenerate fiber
 $W_0=Y_1\cup_E Y_2$ over ${\mathbf 0}\in{\Bbb A}^1$,
 where
  $p_1: Y_1=\Bl_Z X\rightarrow X$ is the blow-up of $X$ along $Z$,
  $p_2: Y_2={\Bbb P}({\cal N}_{Z/X}\oplus{\cal O}_Z)\rightarrow X$
    is a projective space bundle over $Z$, and
  $E={\Bbb P}{\cal N}_{Z/X}$
   with ${\cal N}_{Z/X}$ being the normal bundle of $Z$ in $X$
   is the exceptional divisor of $p_1\,$.

For a fixed $\pi$-ample line bundle ${\cal L}$ on $W$ and
 a topological type $(g,n;d)$,
then the degeneration formula in Sec.~1 applies to this family.
Using the canonical morphism $W/{\Bbb A}^1\rightarrow X$,
 the moduli stack ${\mathfrak M}({\mathfrak W}, (g,n;d))$
 of stable maps from prestable curves into fibers of $W/{\Bbb A}^1$
 constructed in [Li1] is decomposed into a disjoint union of
 collections of irreducible componensts of
 ${\mathfrak M}({\mathfrak W}, (g,n;d))$
 with each collection labelled by an element in 
 $C_{({\cal L},d)}
   := \{\,\beta\in H_2(X;{\Bbb Z})\,:\, {\cal L}\cdot\beta =d\,\}$.
The construction of [Li2] applies to each collection and gives rise
 to a degeneration formula with respect to a curve class
 $\beta\in C_{({\cal L},d)} \subset H_2(X;{\Bbb Z})\,$:
 \begin{eqnarray*}
  \lefteqn{
   \Psi_{(g,n;\,\beta)}^X(\alpha,\zeta) } \\[.6ex]
   && =\; \sum_{\eta\in\overline{\Omega}_{(g,n;\beta)}^{\cal L}}\;
        \frac{{\mathbf m}(\eta)}{|\Eq(\eta)|}\,
          \sum_{j\in K_{\eta}}\,
           \left[ \Psi_{\Gamma_1}^{(Y_1,E)}(j_1^{\ast}\bar{\alpha}(0),
                    \zeta_{\eta,1,j})\,
                  \bullet\,
                  \Psi_{\Gamma_2}^{(Y_2,E)}(j_2^{\ast}\bar{\alpha}(0),
                    \zeta_{\eta,2,j})
           \right]_0\,,
 \end{eqnarray*}
where
 $\bar{\alpha}$ is any flat extension of $\alpha\in A_{\ast}(X)^{\oplus n}$
 to $A_{\ast}(W/{\Bbb A}^1)^{\oplus n}$ and
 $\Omega_{(g,n;\beta)}^{\cal L}\;
  =\; \{\,
       \eta=(\Gamma_1,\Gamma_2,I)
         \in \Omega_{(g,n;d)}^{\cal L}\, |\,
          p_{1\ast}b(\Gamma_1)+p_{2\ast}b(\Gamma_2)=\beta
      \,\}\,$
 is the $\beta$-compatible subset of $\Omega_{(g,n;d)}^{\cal L}\,$.

To remove the possible ${\cal L}$-dependence on the right-hand
 side of the above identity, one can choose ${\cal L}$ on $W$
 to be associated to a sufficiently very ample line bundle on $X$.
For such ${\cal L}$, the set $\Omega_{(g,n;\beta)}^{\cal L}$
 of admissible triples in the above identity is stabilized to
 an ${\cal L}$-independent set $\Omega_{(g,n;\beta)}$ and
the degeneration formula above becomes intrinsic to the topological
 type $(g,n;\beta)$ in the definition of Gromov-Witten invariants
 of $X$.
We now describe $\Omega_{(g,n;\beta)}$, which will be needed in Sec.~3.

Let $\gamma$ denote both the unique curve class
 in $H_2(Y_i;{\Bbb Z})$ that is contracted by $p_i\,$,
 (they corresponds to the same class on $E$).
It has the property that
 $\gamma\cdot E=-1$ on $Y_1$ while $\gamma\cdot E =+1$ on $Y_2$.
A curve class $\beta^{\,\prime}\in H_2(X;{\Bbb Z})$
 has a unique lifting to a curve class $\tilde{\beta}^{\,\prime\,0}$
 in $p_{1\ast}^{\;-1}(\beta^{\,\prime})$ such that all the curve
 classses in $p_{1\ast}^{\;-1}(\beta^{\,\prime})$
 can be written in the form
 $\tilde{\beta}^{\,\prime\,0}+l\gamma$
  for some $l\in {\Bbb Z}_{\ge 0}\,$.
Call $\tilde{\beta}^{\,\prime\,0}$ the minimal lifting of
 $\beta^{\prime}$ to $Y_1$.
When $\beta^{\prime}$ is a curve class in $H_2(Z;{\Bbb Z})$,
 $\beta^{\prime}$ also has a minimal lifting to $Y_2$ with
 the similar property.
Define
$$
 \begin{array}{ccc}
  (H_2(Y_1)\times H_2(Y_2))^0_{\beta}  & :=
   & \parbox[t]{32em}{\raggedright\small
      the set of pairs of curve classes
       $(\tilde{\beta}_1^{\,0},\tilde{\beta}_2^{\,0})
          \in H_2(Y_1;{\Bbb Z})\times H_2(Y_2;{\Bbb Z})$ \newline
      such that
       $p_{i\ast}\tilde{\beta}_i^{\,0} \ne 0$,
       $\tilde{\beta}_i^{\,0}$ is the minimal lifting of 
        $p_{i\ast}\tilde{\beta}_i^{\,0}$ to $Y_i$,   \newline
        $i=1,2\,$, and
       $\;p_{1\ast}\tilde{\beta}_1^{\,0}
          +p_{2\ast}\tilde{\beta}_2^{\,0} =\beta\,$;} \\[7ex]
  (H_2(Y_1))^0_{\beta}  & :=
   & \parbox[t]{32em}{\raggedright\small
      the set of the minimal lifting $\tilde{\beta}^{\,0}$
       of $\beta$ such that
       $E\cdot\tilde{\beta}^{\,0} \ge 0\,$;} \\[2ex]
  (H_2(Y_2))^0_{\beta}  & :=
   & \parbox[t]{32em}{\raggedright\small
      the set of the minimal lifting $\tilde{\beta}^{\,0}$
       of $\beta$ to $Y_2\,$. }
 \end{array}
$$
(By definition, $H_2(Y_2)^0_{\beta}$ is non-empty only
  when $\beta$ is representable by a curve in $Z$.
 Either of $H_2(Y_1)^0_{\beta}$ and $H_2(Y_2)^0_{\beta}$
  is either empty or a singleton.)
With these notations,

{\footnotesize
 $$
  \Omega_{(g,n;\beta)}\;
  =\; \coprod_{\tiny
       (\tilde{\beta}_1^{\,0},\,\tilde{\beta}_2^{\,0})
              \in (H_2(Y_1)\times H_2(Y_2))^0_{\beta}}
         \left\{
          \begin{array}{l}
            \eta=(\Gamma_1,\Gamma_2,I) \\
             \mbox{admissible} \\
             \mbox{triple}     \\
             \mbox{for $Y_1\cup_ E Y_2$}
          \end{array}
          \left|
            \begin{array}{l}
              \bullet\hspace{1ex}
               b(\Gamma_1)=\tilde{\beta}_1^{\,0}+l_1\gamma,\,
               b(\Gamma_2)=\tilde{\beta}_2^{\,0}+l_2\gamma\,, \\
              \hspace{2ex} l_1+l_2 =E\cdot \tilde{\beta}_1^{\,0}\,,\;
                           l_1,\,l_2\in{\footnotesizeBbb Z}_{\ge 0}\,;
                                                                 \\[.6ex]
              \bullet\hspace{1ex}
               g(\eta) = g\,,\;
               n_1+n_2 = n\,; \\[.6ex]
              \bullet\hspace{1ex}
               \sum_{i}\mu_{1,i}\,=\,l_2\,;\\[.6ex]
              \bullet\hspace{1ex}
               I\subset \{1,\,\ldots,\, n\}\,,\; |I|=n_1\,.
            \end{array}
          \right.
         \right\}\;
 $$
 $$
  \hspace{2ex}
  \coprod
  \coprod_{\tiny \tilde{\beta}^{\,0}\in  (H_2(Y_1))^0_{\beta}}
   \left\{
          \begin{array}{l}
            \Gamma_1 :\; \mbox{admissible}   \\
             \mbox{weighted graph}           \\
             \mbox{for $(Y_1,E)$}
          \end{array}
      \left|
            \begin{array}{l}
              \bullet\hspace{1ex}
               b(\Gamma_1)
                = \tilde{\beta}^{\,0}
                  +(E\cdot\tilde{\beta}^{\,0})\gamma\,;\\[.6ex]
              \bullet\hspace{1ex}
               g(\Gamma_1) = g\,,\;
               \mbox{$n$-many legs}\,; \\[.6ex]
              \bullet\hspace{1ex}
               \mbox{no roots}\,.      \\[.6ex]
            \end{array}
      \right.
   \right\}
 $$
 $$
  \hspace{4ex}
  \coprod
  \coprod_{\tiny \tilde{\beta}^{\,0}\in (H_2(Y_2))^0_{\beta}}
   \left\{
          \begin{array}{l}
            \Gamma_2 :\; \mbox{admissible}   \\
             \mbox{weighted graph}           \\
             \mbox{for $(Y_2,E)$}
          \end{array}
      \left|
            \begin{array}{l}
              \bullet\hspace{1ex}
               b(\Gamma_2) = \tilde{\beta}^{\,0}\,;\\[.6ex]
              \bullet\hspace{1ex}
               g(\Gamma_2) = g\,,\;
               \mbox{$n$-many legs}\,; \\[.6ex]
              \bullet\hspace{1ex}
               \mbox{no roots}\,.      \\[.6ex]
            \end{array}
      \right.
   \right\}\,.
 $$
} 

\bigskip

\begin{flushleft}
{\bf Conifold degenerations.}
\end{flushleft}
Conifold degenerations have played some important roles in stringy
 dualities and in understanding how different phases of string theory
 may be connected.
Consider a conifold degeneration given by a projective family
 $\pi:W\rightarrow {\Bbb A}^1$ of $3$-varieties
  with the total space $W$ and fibers $W_t$, $t\ne {\mathbf 0}$,
  smooth and $Y:=W_0$ over ${\mathbf 0}\in {\Bbb A}^1$ a conifold
  with an isolated singularity whose local analytic germ is modelled
  on $\Spec\,{\Bbb C}[[\,x,y,z,w\,]]/(xy-zw)$.
After a base change of degree $2$ and resolution of
 resulting singularities, one obtains a semi-stable reduction
 $\pi^{ss}:W^{ss}\rightarrow {\Bbb A}^1$ of $\pi$
  with $W^{ss}_t=W_{t^2}$, $t\ne {\mathbf 0}$ and
   $p:= \pi^{ss}_0 =p_1\cup_E p_2:
    W^{ss}_0=\widetilde{Y}\cup_E Q \rightarrow Y$,
   where
    $p_1:\widetilde{Y}\rightarrow Y$ is the blow-up of $Y$
     at the conifold singularity,
    $Q$ is the quadric hypersurface in ${\Bbb P}^4$,
    $p_2$ contracts the whole $Q$ to the singulaity of $Y$, and
    $E$ sits in $\widetilde{Y}$ as the exceptional divisor of
     the blow-up $\widetilde{Y}\rightarrow Y$ and in $Q$ from 
     the intersection with a hyperplane of ${\Bbb P}^4\,$.

Fix a $\pi^{ss}$-ample line bundle ${\cal L}$ on $W^{ss}$ and
 a topological type $(g,n;d)$,
then the degeneration formula in Sec.~1 applies to the family
 $\pi^{ss}:W^{ss}\rightarrow {\Bbb A}^1$.
Since $Y$ is topologically a deformation retract of $W_t$
 for any $t\in{\mathbf 0}$, $R^{\bullet}\pi_{\ast}{\Bbb Z}_{W}$
 has a direct summand that is the trivial local system on ${\Bbb A}^1$
 whose fiber is canonically isomorphic to $H_2(W_t; {\Bbb Z})$
 for any $t\in {\Bbb A}^1$.
Using the tautological morphism
 $W^{ss}/{\Bbb A}^1\rightarrow W/{\Bbb A}^1$
 from the semi-stable reduction,
it follows that for any fixed $t_0\ne {\mathbf 0}$,
the stack ${\frak M}({\frak W}^{ss}, (g,n;d))$ of stable morphisms
 from prestable curves of genus $g$, with $n$ marked points, into
 fibers of the universal family of the stack ${\mathfrak W}^{ss}$
 of expanded degenerations associated to $W^{ss}/{\Bbb A}^1$
 with ${\cal L}$-degree $d$ can be decomposed into a disjoint union
 of collections of irreducible componensts of
 ${\mathfrak M}({\mathfrak W}^{ss}, (g,n;d))$
 with each collection labelled by an element of
 $C_{({\cal L},d)}
   := \{\,\beta\in H_2(W_{t_0};{\Bbb Z})\,:\, {\cal L}\cdot\beta =d\,\}$.
(Here we identify $W_{t_0}$ with $W^{ss}_{\sqrt{t_0}}$
 to define ${\cal L}\cdot\beta$.)
The construction of [Li2] applies to each collection and gives rise
 to a degeneration formula with respect to a curve class
 $\beta\in C_{({\cal L},d)} \subset H_2(W_{t_0};{\Bbb Z})\,$:
 \begin{eqnarray*}
  \lefteqn{
   \Psi_{(g,n;\,\beta)}^{W_{t_0}}(\alpha,\zeta) } \\[.6ex]
    &&
    =\; \sum_{\eta\in\overline{\Omega}_{(g,n;\beta)}^{\cal L}}\;
        \frac{{\mathbf m}(\eta)}{|\Eq(\eta)|}\,
         \sum_{j\in K_{\eta}}\,
         \left[ \Psi_{\Gamma_1}^{(\widetilde{Y},E)}
                   (j_1^{\ast}\bar{\alpha}(0),\zeta_{\eta,1,j})\,
                  \bullet\,
                  \Psi_{\Gamma_2}^{(Q,E)}
                   (j_2^{\ast}\bar{\alpha}(0),\zeta_{\eta,2,j})
         \right]_0\,,
 \end{eqnarray*}
where
 $\bar{\alpha}$ is any flat extension of
  $\alpha\in A_{\ast}(W_{t_0})^{\oplus n}
   =A_{\ast}(W^{ss}_{\sqrt{t_0}})^{\oplus n}$
  to $A_{\ast}(W^{ss}/{\Bbb A}^1)^{\oplus n}$ and
 $\Omega_{(g,n;\beta)}^{\cal L}\;
   :=\; \{\,
    \eta=(\Gamma_1,\Gamma_2,I)\in\Omega_{(g,n;d)}^{\cal L}\, |\,
                                   p_{1\ast}b(\Gamma_1)=\beta\,\}\,$
 is the $\beta$-compatible subset of $\Omega_{(g,n;d)}^{\cal L}\,$.

The possible ${\cal L}$-dependence on the right-hand side of the
 above identity can be removed by choosing ${\cal L}$ appropriately
 very ample on $W^{ss}$.
For such ${\cal L}$, the set $\Omega_{(g,n;\beta)}^{\cal L}$
 of admissible triples in the above identity is stabilized to
 an ${\cal L}$-independent set
 $\Omega_{(g,n;\beta)}$, described as follows.
Let
 $\gamma_{2,1}$ and $\gamma_{2,2}$ be the curve classes in
  $H_2(\widetilde{Y};{\Bbb Z})$ represented by the two classes
  of rational curves from the two rulings of $E$ and
 $\gamma_2$ be the curve class that generates
  $H_2(Q;{\Bbb Z})(\simeq {\Bbb Z})$.
(Both $\gamma_{2,1}$ and $\gamma_{2,2}$ become $\gamma_2$
 when passing from $\widetilde{Y}$ to $Q$ via $E$.)
It has the property that
 $\gamma_{2,1}\cdot E=\gamma_{2,2}\cdot E=-1$ on $\widetilde{Y}$
 while $\gamma\cdot E=+1$ on $Q$.
Let $\tilde{\beta}^{\,0}$ be the minimal lifting of $\beta$
 to the set of curve classes in $H_2(\widetilde{Y};{\Bbb Z})$.
It is characterized by the property that any curve classes in
 $p_{1\ast}^{\;-1}(\beta)$ is in the set
 $\tilde{\beta}^{\,0}+{\Bbb Z}_{\ge 0}\gamma_{2,1}
   +{\Bbb Z}_{\ge 0}\gamma_{2,2}\,$.
With these notations,

{\footnotesize
 $$
  \Omega_{(g,n;\beta)}\;
  =\; \left\{
       \begin{array}{l}
         \eta=(\Gamma_1,\Gamma_2,I) \\
          \mbox{admissible} \\
          \mbox{triple}     \\
          \mbox{for $\widetilde{Y}\cup_E Q$}
       \end{array}
       \left|
         \begin{array}{l}
           \bullet\hspace{1ex}
            b(\Gamma_1)
             = \tilde{\beta}^{\,0}
                        +l_{1,1}\gamma_{2,1}+l_{1,2}\gamma_{2,2},\,
            b(\Gamma_2) = l_2\gamma_2\,, \\
           \hspace{2ex} (l_{1,1}+l_{1,2})+l_2
                          =D\cdot \tilde{\beta}^{\,0}\,,\;
                           l_{1,1},\,l_{1,2},\,l_2
                           \in{\footnotesizeBbb Z}_{\ge 0}\,;  \\[.6ex]
           \bullet\hspace{1ex}
            g(\eta) = g\,,\;
            n_1+n_2 = n\,; \\[.6ex]
           \bullet\hspace{1ex}
            \sum_{i}\mu_{1,i}\,=\,l_2\,;\\[.6ex]
           \bullet\hspace{1ex}
            I\subset \{1,\,\ldots,\, n\}\,,\; |I|=n_1\,.
         \end{array}
       \right.
      \right\}\;.
 $$
} 

\noindent
And the degeneration formula above becomes intrinsic to the
 topological type $(g,n;\beta)$ in the definition of
 Gromov-Witten invariants of $W_{t_0}$.

\bigskip

\noindent
{\it Remark 2.2 $[$upshot$\,]$.}
 For our applications, we write the set $\Omega_{(g,n;\,\beta)}$
  of admissible triples $(\Gamma_1,\Gamma_2,I)$ for the degeneration
  associated to blow-ups and to conifold degenerations very explicitly,
  making it look messy.
 However, it should be noted that the essential conditions are only
  $\,p_{1\ast}b(\Gamma_1)+p_{2\ast}b(\Gamma_2)=\beta\,$ and 
  $\,b(\Gamma_1)\cdot E = b(\Gamma_2)\cdot E\,$
  on curve classes.

\bigskip

\section{Transformation of GW-invariants of ${\mathbf 3}$-folds.}

We now employ Sec.~1 and Sec.~2 to study the transformation of
 Gromov-Witten invariants of projective $3$-folds under standard
 flops and small extremal transitions.
The discussion here is in algebro-geometric parallel to that
 of [L-R: Sec.~6] in the symplectic/differential-topological
 category.
Together this gives an algebro-geometric account of
 [Wit5: Sec.~5.5] from the stringy viewpoint and
 [L-R] from the symplectic viewpont.

Recall the following two definitions, which will be used
 in this section:

\bigskip

\noindent
{\bf Definition 3.0.1 [isomorphic power series].} ([L-R].) {\rm
 Two power series $F$ and $G$ are called {\it isomorphic}
 if there exist decompositions $F=F_1+F_2$ and $G=G_1+G_2$
  such that $F_1=G_1$ and that $F_2$ and $G_2$ are related
  to each other by analytic continuations.
}  

\bigskip

\noindent
{\bf Definition 3.0.2 [${\mathbf 3}$-point function].} (cf.\ [C-K].)
{\rm
 Given a smooth projective variety $X$,
 the {\it $3$-point $($correlation$)$ function}
  $$
   \Psi^X\; :\;  A_{\ast}(X)_{\scriptsizeBbb Q}
            \times A_{\ast}(X)_{\scriptsizeBbb Q}
            \times A_{\ast}(X)_{\scriptsizeBbb Q}\;
    \longrightarrow\;
               {\Bbb Q}\,\{\!\{\,H_2(X;{\Bbb Z})\,\}\!\}
  $$
  of cycle classes on $X$ from Gromov-Witten theory
  is defined by
 $$
  \Psi^X(\alpha_1,\alpha_2,\alpha_3)\;
   =\; \sum_{\beta\in H_2(X;{\scriptsizeBbb Z})}\,
        \Psi^X_{(0,3;\,\beta)}(\alpha_1,\,\alpha_2,\,\alpha_3)\,
        q^{\beta}\,.
 $$
} 

\bigskip

\subsection{Gromov-Witten invariants of ${\mathbf 3}$-folds
            under a flop.}

Transformation of Gromov-Witten invariants of projective $3$-folds
 under a standard flop in the algebro-geometric category is explained
 in this subsection.

\bigskip

\begin{flushleft}
{\bf The standard flop of $3$-folds and the associated degenerations.}
\end{flushleft}
Consider the following diagram of a standard flop of
 projective $3$-folds:
 $$
  \begin{array}{cl}
   (\widetilde{X}, E)                             & \\
   \raisebox{1ex}{\scriptsize $p_1$}
     \swarrow\hspace{2em}\searrow
       \raisebox{1ex}{\scriptsize $p_1^{\prime}$}     & \\
   (X,C)\hspace{1.4em}\stackrel{\phi}{\dashrightarrow}\hspace{1em}
   (X^{\prime},C^{\prime})                        &, \\
   \searrow\hspace{2em}\swarrow                   & \\
   \underline{X}                                  &
  \end{array}
 $$
where
 $X$ (resp.\ $X^{\prime}$) is a smooth projective $3$-fold
  with an embedded smooth curve
  $C$ (resp.\ $C^{\prime}$) $\simeq {\Bbb P}^1$
  whose normal sheaf ${\cal N}_{C/X}$
  (resp.\ ${\cal N}_{C^{\prime}/X^{\prime}}$)
  is isomorphic to
  ${\cal O}_{{\scriptsizeBbb P}^1}(-1)
               \oplus{\cal O}_{{\scriptsizeBbb P}^1}(-1)$;
 $(\tilde{X},E)$ is the common blow-up of $X$ along $C$
   and of $X^{\prime}$ along $C^{\prime}$ with the exceptional
   divisor $E\simeq {\Bbb P}^1\times{\Bbb P}^1$;
 $X\rightarrow \underline{X}$
  (resp.\ $X^{\prime}\rightarrow \underline{X}$)
  is the morphism that contracts exactly $C$ (resp.\ $C^{\prime}$);
 the birational map $\phi: X\dashrightarrow X^{\prime}$ is
  the composition of blowing up $X$ along $C$ and then blowing down
  the ruling of $E$ not contracted by $\widetilde{X}\rightarrow X$,
 $\phi$ (or simply $X^{\prime}$ when the diagram is understood
  implicitly) is called the {\it flop} of $X\rightarrow \underline{X}$
  (or simply $X$).

There are two degenerations associated to the above diagram:
 $$
  \pi:W:=\Bl_{\,C\times{\mathbf 0}}\,X\times{\Bbb A}^1
   \rightarrow {\Bbb A}^1 
  \hspace{1em}\mbox{and}\hspace{1em}
  \pi^{\prime}:
  W^{\prime}
  :=\Bl_{\,C^{\prime}\times{\mathbf 0}}\, X^{\prime}\times {\Bbb A}^1
                                           \rightarrow {\Bbb A}^1\,.
 $$
Both $\pi$ and $\pi^{\prime}$ are constant families
 except over ${\mathbf 0}$, where the fiber becomes
  $$
   W_0\; 
   =\; \widetilde{X}\cup_E {\Bbb P}({\cal N}_{C/X}\oplus {\cal O}_C)\;
   =: \widetilde{X}\cup_E Y \hspace{1em}
  \mbox{for $\pi$}
  $$
 and
  $$
   W^{\prime}_0\;
   =\; \widetilde{X}\cup_E {\Bbb P}({\cal N}_{C^{\prime}/X^{\prime}}
                                     \oplus {\cal O}_{C^{\prime}})\;
   =: \widetilde{X}\cup_E Y^{\prime}  \hspace{1em}
  \mbox{for $\pi^{\prime}$}\,.
  $$
 Here $Y$ and $Y^{\prime}$ are isomorphic but the gluing to
  $\widetilde{X}$ along $E$ differs by an exchange of the two
  ${\Bbb P}^1$ factors in $E\simeq {\Bbb P}^1\times{\Bbb P}^1$.
 Such an automorphism on $E$ as a subvariety in $Y$ or $Y^{\prime}$
  does not extend to an automorphism on the whole $Y$ or $Y^{\prime}$;
 thus in general $W_0$ and $W^{\prime}_0$ are not isomorphic.

The birational map $\phi:X\rightarrow X^{\prime}$ induces
 degree-preserving homomorphisms on Chow groups
 $$
  \phi_{\ast}:A_{\ast}(X)\rightarrow A_{\ast}(X^{\prime})
   \hspace{1em}\mbox{and}\hspace{1em}
  \phi^{-1}_{\ast}:
    A_{\ast}(X^{\prime})\rightarrow A_{\ast}(X)
 $$
 via intersection with the (closed) graph $\Gamma_{\phi}$ of
  $\phi$ in $X\times X^{\prime}$.
Since $\phi$ is a curve-surgery on a $3$-fold,
 $\phi_{\ast}: A_3(X)\rightarrow A_3(X^{\prime}),
  A_2(X)\rightarrow A_2(X^{\prime})$
 are isomorphisms.
For a curve class $\gamma\in A_1(X)$ with $\gamma\ne [C]$,
 using Chow's moving lemma, one can represent $\gamma$ by a linear
 combination of $1$-cycles that are disjoint from $C$.
This gives a representation of $\phi_{\ast}(\gamma)$ by a linear
 combination of $1$-cycles disjoint with $C^{\prime}$.
Since $X^{\prime}$ is projective, $C^{\prime}$ is the only
 locus contracted by $X^{\prime}\rightarrow \underline{X}$,
 and the preperty of being contracted is a numerical property,
 one concludes that $\phi_{\ast}([C])\ne [C^{\prime}]$.
The same reasoning implies that $\phi_{\ast}([C]) = k[C^{\prime}]$
 for some $k\in {\Bbb Z}$.
To compute $k$, observe that,
 though $\Gamma_{\phi}$ and $C\times X^{\prime}$
  have excess intersection along $E$ of dimension $1$ larger
  than expected in $X\times X^{\prime}$,
 both the embeddings
  $\Gamma_{\phi}\cap (C\times X^{\prime})\simeq E
   \hookrightarrow X\times X^{\prime}$ and
  $C\times X^{\prime}\hookrightarrow X\times X^{\prime}$
  are regular.
The excess intersection formula gives then
 $\Gamma_{\phi}\cdot(C\times X^{\prime})
  = c_1({\cal N}_{(C\times X^{\prime})/(X\times X^{\prime})}
                 / {\cal N}_{E/\tilde{X}})\cap [E]$,
 which is the class $[C]-[C^{\prime}]$
 on $E\simeq C\times C^{\prime}$.
This class projects to $-[C^{\prime}]$ on $X^{\prime}$.
Thus, $\phi_{\ast}([C])=-[C^{\prime}]$.
This argument shows also that
 $\phi_{\ast}: A_1(X)\rightarrow A_1(X^{\prime})$
 is an isomorphism.
Finally, one has a homomorphism
 $\phi_{\ast}:A_0(X)\rightarrow A_0(X^{\prime})$.

For the degenerate fiber
 $W_0=\widetilde{X}\cup_E Y$ in $W/{\Bbb A}^1$,
 the Chow rings $A_{\ast}\widetilde{X}$ are determined by
 $A_{\ast}(X)$ via blow-up formulas.
The Chow ring $A_{\ast}(Y)$ is generated by $w:=[E]$ and
 the fiber class $v:=[{\Bbb P}^2]$ of $Y\rightarrow C$.
Explicitly, $A_{\ast}(Y)\simeq {\Bbb Z}[v,w]/(v^2, w^3-2vw^2)$
 by either the projective space bundle or the toric computation.
A class $[\zeta_1,\zeta_2]$ in $A_{\ast}(W_0)$
 (with $\zeta_1$ a linear combination of cycles in $\widetilde{X}$
  and $\zeta_2$ a linear combination of cycles in $Y$)
 that comes from a relative class
  $\alpha\in H^0(R^{\bullet}\pi_{\ast}{\Bbb Q}_W)$
 must satisfy the predeformability condition
  $\zeta_1\cdot E = \zeta_2\cdot E\,$.

Consider now the classes in $H^0(R^{\bullet}\pi_{\ast}{\Bbb Q}_W)$
 that are representable by linear combinations of relative cycles
 on $W/{\Bbb A}^1$.
We will call these classes algebraic.
For an algebraic class $\alpha$ in $H^0(R^6\pi_{\ast}{\Bbb Q}_W)$
 and $H^0(R^4\pi_{\ast}{\Bbb Q}_W)$,
since $W-W_0$ is a trivial family with fiber $X$ over
 ${\Bbb A}^1-{\mathbf 0}$,
$\alpha$ determines a unique class $\alpha_t := \alpha|_{W_t}$
 ($t\ne {\mathbf 0}$) in $A_{\ast}(X)$.
Using Chow's moving lemma and taking completion,
 $\alpha$ determines a class
 $\alpha^{\prime}\in H^0(R^4\pi_{\ast}{\Bbb Q}_W)$
 that coincides with $\alpha$ over ${\Bbb A}^1-{\mathbf 0}$
 and is represented by a linear combination of relative cycles that
  intersect $W_0$ only at the $\widetilde{X}$ component.
For an algebraic class $\alpha\in H^0(R^2\pi_{\ast}{\Bbb Q}_W)$,
 same argument concludes that $\alpha$ determines a class
 $\alpha^{\prime}\in H^0(R^2\pi_{\ast}{\Bbb Q}_W)$ that 
 coincides with $\alpha$ over ${\Bbb A}^1-{\mathbf 0}$
 and can be represented by
 ${p_1^{-1}}_{\ast}(\alpha_t) + (\alpha_t\cdot C)[{\Bbb P}^2]$
 in $A^1(\widetilde{X}\cup_E Y)$.
Finally, $H^0(R^0\pi_{\ast}{\Bbb Q}_W)$ is generated by $[W]$.

When one is restricted to the $\widetilde{X}$ part of cycles,
 the above homomorphism coincides with
 ${p_1^{-1}}_{\ast}: A_{\ast}(X)\rightarrow A_{\ast}(\widetilde{X})$.
Together with a similar construction for the degeneration
 $\pi^{\prime}:W^{\prime}\rightarrow {\Bbb A}^1$,
 one has a commutative diagram:
$$
 \begin{array}{cl}
  A_{\ast}(\widetilde{X})                             & \\
  \raisebox{1ex}{\scriptsize ${p_1^{-1}}_{\ast}$}
    \nearrow\hspace{2em}\nwarrow
      \raisebox{1ex}{\scriptsize ${{p_1^{\prime}}^{-1}}_{\ast}$}   & \\
  A_{\ast}(X)
   \hspace{1.4em}\stackrel{\phi_{\ast}}{\longrightarrow}\hspace{1em}
  A_{\ast}(X^{\prime})                        &.
 \end{array}
$$
There is a cycle class map
 $A_{\ast}(\,\cdot\,)\rightarrow H_{\ast}(\,\cdot\,)$
and the above diegram is compatibility with the topological results
on homomorphisms among the associated $H_{\ast}(\,\cdot\,)$ in [L-R].

\bigskip

\begin{flushleft}
{\bf GW-invariants of $X$, $X^{\prime}$ in terms of
     relative GW-invariants of $(\widetilde{X},E)$, $(Y,E)$.}
\end{flushleft}
Recall
 the moduli stack ${\mathfrak M}({\mathfrak W}, (g,n;\beta))$
  of stable maps of topological type $(g,n;\beta)$ to fibers of
  the universal family of the stack ${\mathfrak W}$ of expanded
  degenerations associated to $W/{\Bbb A}^1$
 and the set $\Omega_{(g,n;\beta)}$ of admissible triples that
  appears in the degeneration formula and is intrinsic to
  $(g,n;\beta)$.

For a general projective degeneration $W\rightarrow {\Bbb A}^1$
 with $W_0=Y_1\cup_D Y_2$ discuessed in [Li1] and
with the extra condition that $H_2(W_t;{\Bbb Z})$
 is canonically isomorphic to $H_2(W_{t_0};{\Bbb Z})$
 for some $t_0\ne {\mathbf 0}$,
similar discussions as in [L-Y1] always give a refined
 degeneration formula from [Li2] for curve classes
 $\beta\in H_2(W_{t_0};{\Bbb Z})$.
In such cases, suppose that $\beta$ is decomposed to
 $\beta_1+\beta_2$ in $H_2(W_0;{\Bbb Z})$
 (when represented in $H_2(Y_1;{\Bbb Z})\oplus H_2(Y_2;{\Bbb Z})$,
   $(\beta_1,\beta_2)$ may not be unique,
   simply choose any representative) and
 that $(\Gamma_1,\Gamma_2,I)$ is an admissible triple in
  $\Omega_{(g,n;\beta)}$ such that
   $b(\Gamma_i)=\beta_i$ and that 
   $\Gamma_i$ has $n_i$ legs and $r$ roots, $i=1,\,2\,$. 
Let $g_i:= g(\Gamma_i)$.
Then
 the following additivity relation from the (complex)
  dimension of the virtual fundamental classes
  $[\overline{\cal M}_{g,n}(W_{t_0},\beta)]^{\virt}$ and 
  $[{\mathfrak M}((Y_i,D;\Gamma_i))]^{\virt}\,$, $i=1,\,2\,$, holds:
 \begin{eqnarray*}
  \lefteqn{(1-g)(\dimm W_{t_0} -3) - \beta\cdot K_{W_{t_0}} + n} \\
   && =\;  (1-g_1)(\dimm Y_1 -3) - \beta_1\cdot K_{Y_1}
             + n_1 + (r-\beta_1\cdot D)  \\
   && \hspace{2em}
          + \, (1-g_2)(\dimm Y_2 -3) - \beta_2\cdot K_{Y_2}
             + n_2 + (r-\beta_2\cdot D)
          - r (\dimm D)\,.
 \end{eqnarray*}
By construction,
 $\dimm W_t=\dimm Y_1=\dimm Y_2 = \dimm D +1$,
 $g = g_1+g_2 + r+1
      -|V(\Gamma_1\,\raisebox{.2ex}{\scriptsize $\coprod$}\, \Gamma_2)|$,
 $n=n_1+n_2$, 
 $\beta_1\cdot D= \beta_2\cdot D \ge 0$, and
 $0\le r \le \beta_1\cdot D\,$.

For the degeneration $W/{\Bbb A}^1$ with $W_0=\widetilde{X}\cup_E Y$
 from the standard flop of $3$-folds,
 $p_1:\widetilde{X}\rightarrow X$ is the blow-up of $X$
  along $C$ with exceptional divisor $E$,
  hence $K_{\widetilde{X}}=p_1^{\ast}K_X+E$
  by computing the Jacobian of $p_1$ at a general point of $E$;
 $Y$ is a toric $3$-fold, whose fan structure gives $K_Y=-3E$.
The additivity relation thus simplifies to
$$
 -\beta\cdot K_X + n\;
 =\;   \left(\, -\beta_1\cdot p_1^{\ast}K_X + n_1 + r - 2\beta_1\cdot E\,
       \right)\,
   +\, \left(\, 2\beta_2\cdot E + n_2 + r \, \right)\,
   -\, 2r\,.
$$
There are two cases.

\bigskip

\noindent
{\it Case $(a)\,$}:
{\it the virtual dimension
     $\vdim\overline{\cal M}_{g,\,0}(X,\beta)=0$.}
 In this case, one has the identity
   $$
    0\;
    =\;   \left(\, -\beta_1\cdot p_1^{\ast}K_X + r - 2\beta_1\cdot E\,
          \right)\,
      +\, \left(\, 2\beta_2\cdot E + r \right)\,
      -\, 2r\,,
   $$
  which implies that $\beta_1\cdot p_1^{\ast}K_X=0$.
  Thus,
   $$
    \vdim\, {\mathfrak M}(\widetilde{X},E;\Gamma_1)\;
     =\; r-2\beta_1\cdot E
     \hspace{1em}\mbox{and}\hspace{1em}
    \vdim\, {\mathfrak M}(Y,E;\Gamma_2)\; =\; 2\beta_2\cdot E + r\,.
   $$

For $\beta$ not a multiple of $[C]$, it always holds that
 $\beta_1\ne 0$.
In this case, the only situation that
 $\vdim\overline{\cal M}_{g_1,0}(\widetilde{X},E;\beta_1,r)\ge 0$
 is when $\beta_1\cdot E=r=0$.
This implies that $\beta_2=0$.
It follows from Sec.~2 that in this case
 elements in $\Omega_{(g,0;\beta)}$ correspond to a unique curve
  class $\tilde{\beta}$
  in $H_2(\widetilde{X};{\Bbb Z})$, characterized by the conditions
  $p_{1\ast}\tilde{\beta}=\beta$ and $\tilde{\beta}\cdot E=0$.
This class coincides with the class associated to $\beta$ from
 Chow's moving lemma.
Thus, for $\beta$ not a multiple of $[C]$, the Gromov-Witten
   invariant $\Psi^X_{(g,0;\,\beta)}$ of $X$ is equal to
   the relative Gromov-Witten invariant
   $\Psi^{(\widetilde{X},E)}_{(g,0;\,\tilde{\beta})}$
   of $(\widetilde{X},E)$.

If $\beta$ is a positive multiple of $[C]$,
 then if $\beta_1\ne 0$, it must lie in $E$ with $\beta_1\cdot E<0$.
Such $\beta_1$ does not occur in any admissible triple
 in $\Omega_{g,0;\beta}$.
It follows from Sec.~2 that in this case
 elements in $\Omega_{(g,0;\beta)}$ correspond to
 a unique curve class $\tilde{\beta}$ in $H_2(Y;{\Bbb Z})$,
 characterized by the conditions
 $p_{2\ast}\tilde{\beta}=\beta$ and $\tilde{\beta}\cdot E=0$.
Recall that $Y$ is a toric variety.
Such $\tilde{\beta}$ on $Y$ is represented by a multiple of
 the unique toric invariant curve $\overline{C}:={\Bbb P}^1$
 not intersecting $E$ and with normal bundle
 ${\cal O}_{{\scriptsizeBbb P}^1}(-1)
  \oplus {\cal O}_{{\scriptsizeBbb P}^1}(-1)$.
The related Gromov-Witten invariant $\Psi^X_{(g,0;\,\beta)}$ of $X$
 is identical to the absolute Gromov-Witten invariant
 $\Psi^Y_{(g,0;\,\tilde{\beta})}$.
  
\bigskip

\noindent
{\it Case $(b)\,$}:
{\it $\vdim\overline{\cal M}_{g,0}(X,\beta)>0$.}
In this case, one considers the Gromov-Witten invariants
 $\Psi^X_{(g,n;\,\beta)}(\alpha_1,\,\cdots,\,\alpha_n)$
 with the summation of the codimension of the cycle $\alpha_i$ on
 $X$ equal to $\vdim\overline{\cal M}_{g,n}(X,\beta)$.
Since our purpose is to compare Gromov-Witten invariants
 under a standard flop, by the Fundamental Axiom and
 the Divisor Class Axiom,
one only needs to consider $\alpha_i\in A_0(X)$ or $A_1(X)$.
From the earlier discussion in this section, each such class
 determines a class in
 $H^0(R^{\,6}\pi_{\ast}{\Bbb Q}_W)+H^0(R^4\pi_{\ast}{\Bbb Q}_W)$
 representable by a relative cycle on $W/{\Bbb A}^1$
 that intersects $W_0$ only at $\widetilde{X}$.
The only contribution to
 $\Psi^X_{(g,n;\,\beta)}(\alpha_1,\,\cdots,\,\alpha_n)$
 from the summands of the degeneration formula comes from
 those labelled by admissible triples with $\beta_2=0$.
As in Case (a), when $\beta$ is not a multiple of $[C]$,
 there is a unique such admissible triple
  $(\Gamma_1, \emptyset, \emptyset)$ in $\Omega_{(g,n;\beta)}$
 with $b(\Gamma_1)=\tilde{\beta}$.
This $\tilde{\beta}$ is again characterized by
 $p_{1\ast}\tilde{\beta}=\beta$ and $\tilde{\beta}\cdot E=0\,$.
Consequently,
 $\Psi^X_{(g,n;\,\beta)}(\alpha_1,\,\cdots,\,\alpha_n)
  = \Psi^{(\widetilde{X},E)}_{(g,n;\,\tilde{\beta})}
     ({p_1^{-1}}_{\ast}\alpha_1,\,\cdots,\, {p_1^{-1}}_{\ast}\alpha_n)$.
If $\beta$ is a multiple of $[C]$, there is no such admissible triple
 and $\Psi^X_{(g,n;\,\beta)}(\alpha_1,\,\cdots,\,\alpha_n)=0$,
 which is consistent with the vanishing for the dimensional reason
 in this situation.

\bigskip

The discussion for $X^{\prime}$ is identical.

\bigskip

\begin{flushleft}
{\bf Gromov-Witten invariants under the standard flop.}
\end{flushleft}
Comparing Gromov-Witten invariants of $X$ and $X^{\prime}$ both
 to the relative Gromov-Witten invariants of $(\widetilde{X},E)$
 or the absolute Gromov-Witten invariants of $Y$
 from the discussion above,
recalling the commutative diagram:
 $$
  \begin{array}{cl}
   A_{\ast}(\widetilde{X})                             & \\
   \raisebox{1ex}{\scriptsize ${p_1^{-1}}_{\ast}$}
     \nearrow\hspace{2em}\nwarrow
       \raisebox{1ex}{\scriptsize ${{p_1^{\prime}}^{-1}}_{\ast}$}   & \\
   A_{\ast}(X)
    \hspace{1.4em}\stackrel{\phi_{\ast}}{\longrightarrow}\hspace{1em}
   A_{\ast}(X^{\prime})                        &,
  \end{array}
 $$
and together with the axioms of Gromov-Witten invariants and
Chow's moving lemma gives the following relations:

\bigskip

\noindent
{\bf Theorem 3.1.1 [GW-invariant under flop].}
([L-R: Theorem A and Corollary A.1].)
{\it
 Let $\beta\in H_2(X;{\Bbb Z})$ and
  $\alpha_1\,,\,\cdots\,,\,\alpha_n
                  \in A_{\ast}(X)_{\scriptsizeBbb Q}$.
 If $\beta$ is not a multiple of $[C]$, then
  $$
   \Psi^X_{(g,n;\,\beta)}(\alpha_1,\,\cdots,\,\alpha_n)
   =\Psi^{X^{\prime}}_{(g,n;\,\phi_{\ast}(\beta))}
     (\phi_{\ast}\alpha_1,\,\cdots,\,\phi_{\ast}\alpha_n)
  $$
 while
  $\Psi^X_{(g,0;\,m[C])}=\Psi^{X^{\prime}}_{(g,0;\,m[C^{\prime}])}$.
 If $X$ is a Calabi-Yau $3$-fold,
  then the first identity is reduced to
  $\Psi^X_{(g,0;\,\beta)}
         =\Psi^{X^{\prime}}_{(g,0;\,\phi_{\ast}(\beta))}$.
} 

\bigskip

This theorem implies the transformation law of $3$-point
 functions under a standard flop.
The discussion below in rational cycle classes is similar to
 [L-R] in differential forms and their Poincar\'{e} dual.

The $3$-point function can be decomposed into three parts:
\begin{eqnarray*}
 \lefteqn{\Psi^X(\alpha_1,\alpha_2,\alpha_3)} \\
  &&
  =\; (\alpha_1\cdot\alpha_2\cdot\alpha_3)_0
      + \sum_{\beta\ne 0, m[C]}\,
         \Psi^X_{(0,3;\,\beta)}(\alpha_1,\alpha_2,\alpha_3)\,q^{\beta}
      + \sum_{\beta=m[C], m>0}
         \Psi^X_{(0,3;\,\beta)}(\alpha_1,\alpha_2,\alpha_3)\,q^{\beta}\,,
\end{eqnarray*}
where the first term takes only the $0$-dimensional component of
 the triple intersection product.
The power series is invariant under permutations of $\alpha_i$'s
 and we will assume that
 $\codimm\alpha_1\ge\codimm\alpha_2\ge \codimm\alpha_3$.
Theorem 3.1.1 implies that the second term is invariant under
 $q^{\beta}\rightarrow q^{\,\phi_{\ast}(\beta)}$.

For dimensional reason, the third term vanishes except for
 $\alpha_i\in A^1(M)$.
When there are no insertions,
 $\Psi^X_{(0,0;\,m[C])}=\Psi^Y_{(0,0;\,m[\overline{C}])}=1/m^3$
  by a localization computation,
(see [C-K] for a review and references therein on this result
  by several independent authors;
 in [L-L-Y] this is a by-product of the Mirror Principle framework,
 in which a linearized moduli space for stable maps from [Wit5]
  is integrated into both the Gromov-Witten theory and
  the localization computations).
Consequently, the third term is given by 
 $$
  \sum_{\beta=m[C], m>0}
   \Psi^X_{(0,3;\,\beta)}(\alpha_1,\alpha_2,\alpha_3)\,q^{\beta}\;
  =\;([C]\cdot\alpha_1)\, ([C]\cdot\alpha_2)\, ([C]\cdot\alpha_3)\;
      \frac{ q^{[C]} }{ 1-q^{[C]} }\,.
 $$
The first term vanishes except for
 $(\codimm\alpha_1,\codimm\alpha_2,\codimm\alpha_3)=(1,1,1)$,
  $(2,1,0)$, or $(3,0,0)$.
Compare now the first and the third term on the $X$ and
 the $X^{\prime}$ side.

For $(\codimm\alpha_1,\codimm\alpha_2,\codimm\alpha_3)=(3,0,0)$
 or $(2,1,0)$, the third term vanishes and
 the identity
  $\alpha_1\cdot\alpha_2\cdot\alpha_3
   =\phi_{\ast}(\alpha_1)\cdot
    \phi_{\ast}(\alpha_2)\cdot\phi_{\ast}(\alpha_3)$
 follows from Chow's moving lemma, which renders the intersections
 a same combination of $0$-cycles in
 $X-C\stackrel{\phi}{\simeq} X^{\prime}-C^{\prime}$.

For $(\codimm\alpha_1,\codimm\alpha_2,\codimm\alpha_3)=(1,1,1)$,
usineg Chow's moving lemma, one can assume that
 the two representing cycles for each of the following pairs meet
 properly
 $$
  \begin{array}{c}
   (\,\alpha_1\,,\, [C]\,)\,, \hspace{2em}
   (\,\alpha_2\,,\, \alpha_1+[C]+ [\alpha_1\cap C]\,)\,, \\[1ex]
  (\,\alpha_3\,,\,
     \alpha_1+\alpha_2+[C]+[\alpha_1\cap\alpha_2]+
           [\alpha_1\cap C] + [\alpha_2\cap C]\,)\,.
  \end{array}
 $$
(For simplicity of notation, we identify $\alpha_i$ with
 their representing cycles in the discussion below.)
This implies that $\alpha_1\cdot\alpha_2\cdot\alpha_3$
 is represented by a linear combination of $0$-cycles disjoint
 from $C$ and that 
$\phi_{\ast}(\alpha_1)$, $\phi_{\ast}(\alpha_2)$,
 $\phi_{\ast}(\alpha_3)$ meet properly with each other
 except along $C^{\prime}$, where an excess intersection of
 the three may occur.
Under this choice of representing cycles,
 every two of $\phi_{\ast}(\alpha)$, $\phi_{\ast}(\alpha_2)$,
  and $\phi_{\ast}(\alpha_3)$ meet properly along $C^{\prime}$
  in $X^{\prime}$ if $C^{\prime}$ contains their intersection points.
This implies that 
 \begin{eqnarray*}
  \lefteqn{\phi_{\ast}(\alpha_1)\cdot\phi_{\ast}(\alpha_2)
   \cdot\phi_{\ast}(\alpha_3)}\\
  && =\; \alpha_1\cdot\alpha_2\cdot\alpha_3\,
          +\, (\,\mbox{contribution from the excess intersection
                       along $C^{\prime}$}\,)
 \end{eqnarray*}
 for the triple intersecion.
The intersection multiplicity
  $i(C^{\prime},
      \phi_{\ast}(\alpha_1)\cdot \phi_{\ast}(\alpha_2);
      X^{\prime})$
 of $\phi_{\ast}(\alpha_1)$ and $\phi_{\ast}(\alpha_2)$
 along $C^{\prime}$ in $X^{\prime}$ is given by
  $(-[C]\cdot\alpha_1)(-[C]\cdot\alpha_2)$.
Since $[C^{\prime}]\cdot \phi_{\ast}(\alpha_3)= -[C]\cdot\alpha_3$,
 one concludes that 
 $$
  \phi_{\ast}(\alpha_1)\cdot\phi_{\ast}(\alpha_2)
   \cdot\phi_{\ast}(\alpha_3)\;
   =\; \alpha_1\cdot\alpha_2\cdot\alpha_3\,
      -\, ([C]\cdot\alpha_1)\,
          ([C]\cdot\alpha_2)\, ([C]\cdot\alpha_3)\,.
 $$
Finally, note that
 $\phi_{\ast}(q^{[C]}/(1-q^{[C]}))
   = q^{\phi_{\ast}[C]}/(1-q^{\phi_{\ast}[C]})
   = q^{-[C^{\prime}]}/(1-q^{-[C^{\prime}]})$
 is isomorphic to $-1 - q^{[C^{\prime}]}/(1-q^{[C^{\prime}]})$
 by an analytic continuation.
Combining all these and re-writing the triple intersection above as
 \begin{eqnarray*}
  \lefteqn{
    \phi_{\ast}(\alpha_1\cdot\alpha_2\cdot\alpha_3)  } \\[.6ex]
    && =\; \phi_{\ast}(\alpha_1)\cdot\phi_{\ast}(\alpha_2)
             \cdot\phi_{\ast}(\alpha_3)\,
          -\, ([C^{\prime}]\cdot\phi_{\ast}(\alpha_1))\,
              ([C^{\prime}]\cdot\phi_{\ast}(\alpha_2))\,
              ([C^{\prime}]\cdot\phi_{\ast}(\alpha_3))\,,
 \end{eqnarray*}
 one concludes that 
\begin{eqnarray*}
 \lefteqn{
  \phi_{\ast}\left(\,\alpha_1\cdot\alpha_2\cdot\alpha_3
   + ([C]\cdot\alpha_1)([C]\cdot\alpha_2)([C]\cdot\alpha_3)\,
     \frac{q^{[C]}}{1-q^{[C]}}\,
               \right)            }\\
 && \stackrel{a.c.}{=}\;
         \phi_{\ast}(\alpha_1)\cdot\phi_{\ast}(\alpha_2)
                                \cdot\phi_{\ast}(\alpha_3)
        + ([C^{\prime}]\cdot\phi_{\ast}(\alpha_1))
          ([C^{\prime}]\cdot\phi_{\ast}(\alpha_2))
          ([C^{\prime}]\cdot\phi_{\ast}(\alpha_3))\,
          \frac{q^{[C^{\prime}]}}{1-q^{[C^{\prime}]}}\,,
\end{eqnarray*}
where {\it a.c.} stands for ``analytic continuation".
The second line of this identity is exactly
 the summation of the first and the third part of the similar
 decomposition of the $3$-point function
 $\Psi^{X^{\prime}}(\phi_{\ast}\alpha_1, \phi_{\ast}\alpha_2,
                    \phi_{\ast}\alpha_3)$.
Recall Definition 3.0.1 of isomorphic power series.
The whole discussions thus imply that:

\bigskip

\noindent
{\bf Corollary 3.1.2 [${\mathbf 3}$-point function under flop].}
 ([L-R: Corollary A.2].) {\it
 The $3$-point function
  $\Psi^{X^{\prime}}(\,\phi_{\ast}(\alpha_1),
                       \phi_{\ast}(\alpha_2),
                       \phi_{\ast}(\alpha_3)\,)$ of $X^{\prime}$
 is isomorphic to  
 $\phi_{\ast}(\Psi^{X}(\alpha_1,\alpha_2,\alpha_3))\,$.
} 

\bigskip

\noindent
{\bf Corollary 3.1.3 $[$Witten's wall-crossing formula$\,]$.}
Witten's wall-crossing formula [Wit5: Eq.~(5.48)] of Yukawa couplings
 of the massless supermultiplets associated to $H^{1,1}(\,\cdot\,)$
 in the $4$-dimensional low-energy effective field theory from
 compactification of a superstring model can be completed to
 isomorphisms of 3-point functions of Gromov-Witten
 invariants for $3$-folds related by a standard flop
 in [L-R: Corollary A.2], cf.\ Corrollary 3.1.2.
See Sec.~A.1 for more explanations on this stringy aspect of
 the problem.

\bigskip

\subsection{Gromov-Witten invariants of ${\mathbf 3}$-folds
            under a small extremal transition.}

In the diagram of a standard flop of $3$-folds at the beginning
 of Sec.~3.1,
the singular variety $\underline{X}$ from contracting a curve
 $C\simeq {\Bbb P}^1$ in $X$ with normal sheaf
 ${\cal N}_{C/X}
   \simeq {\cal O}_{{\scriptsizeBbb P}^1}(-1)
          \oplus {\cal O}_{{\scriptsizeBbb P}^1}(-1)$
 is a conifold.
We will assume that $\underline{X}$ arises also from a degeneration
 $\pi^{\prime\prime}:W^{\prime\prime}\rightarrow {\Bbb A}^1$
 with $W^{\prime\prime}$ smooth, $\pi^{\prime\prime}$ projective,
 and $W^{\prime\prime}_0=\underline{X}$ the fiber over
 ${\mathbf 0}\in{\Bbb A}^1$.
In this subsection, we explain how Gromov-Witten invariants of
 $X^{\prime\prime}:= W^{\prime\prime}_{t_0}$,
 for a fixed $t_0\ne{\mathbf 0}$, are related to the Gromov-Witten
 invariants of $X$ in the algebro-geometric setting.

\bigskip

\begin{flushleft}
{\bf Correspondences of algebraic classes.}
\end{flushleft}
Recall from Sec.~2 the semi-stable reduction
 $\pi^{ss}:W^{ss}\rightarrow {\Bbb A}^1$ of $\pi^{\prime\prime}$
 from a degree-$2$ base change with the fiber over ${\mathbf 0}$
 the gluing $W^{ss}_0=\widetilde{X}\cup _E Q$,
  where $Q$ is the smooth quadric hypersurface in ${\Bbb P}^4$
   with $E\subset Q$ from a hyperplane section.
The monodromy action of $\pi_1({\Bbb A}^1-{\mathbf 0})$ is trivial
 on $H_k(X^{\prime\prime};{\Bbb Q})$ except for $k=3$,
  which is not relevant in the current algebro-geometric setting.
After the semi-stable reduction, the new monodromy action remains
 trivial on
 $H_k(W^{ss}_{\sqrt{t_0}};{\Bbb Q})
  =H_k(X^{\prime\prime};{\Bbb Q})$ for $k\ne 3$.

Given a class
 $\bar{\alpha}\in A_{\ast}(W^{\prime\prime})$
 reprented by a relative cycle (still denoted by $\bar{\alpha}$)
 of relative dimension $k\le 2$, using Chow's moving lemma, one may
 assume that $\bar{\alpha}$ does not contain the conifold singularity
 in the fiber $W^{\prime\prime}_0=\underline{X}$.
Then the lifting $\bar{\alpha}^{ss}$ of $\bar{\alpha}$ to $W^{ss}$
 induces a class
 $\alpha:=\bar{\alpha}_{t_0}\in A_k(X^{\prime\prime})$,
 and a class
 $\tilde{\alpha}:=\bar{\alpha}_0\in A_k(\widetilde{X})$
 by intersecting $\bar{\alpha}^{ss}$ with $W^{ss}_{\sqrt{t_0}}$ and
 $W^{ss}_0=\widetilde{X}\cup_E Q$ respectively.
(The latter intersection is disjoint from $Q$ by our choice of
 the representative of $\bar{\alpha}$.)
In this way, the subgroup of
 $H^0(R^{\,6-2k}\pi^{\prime\prime}_{\ast}\,
                        {\Bbb Q}_{W^{\prime\prime}})$
 that consists of elements representable by relative cycles on
 $W^{\prime\prime}/{\Bbb A}^1$ induces a correspondence
 $$
  \psi_{\ast}\;:\; A_k(X^{\prime\prime}) \;
                   \vdash\; A_k(\widetilde{X})
   \hspace{2em}\mbox{with}\hspace{2em}
  \alpha\;\leadsto\; \tilde{\alpha}\,.
 $$
(Here, unlike in the smooth category, there exists no rational map
 $\psi:X^{\prime\prime}\dashrightarrow \widetilde{X}$.
 The lower $\ast$ is added only for a notation syncronization.)
Define $\psi_{\ast}:A_3(X^{\prime\prime})\vdash A_3(\widetilde{X})$
  directly by $[X^{\prime\prime}]\mapsto [\widetilde{X}]$.
Let $A_{\ast}(X^{\prime\prime})^{\circ}$ be the subgroup of
 $A_{\ast}(X^{\prime\prime})$ that consists of elements
  whose image under $\psi_{\ast}$ is non-empty.
Recall the blow-up $p_1:\widetilde{X}\rightarrow X$ from Sec.~3.1.
Then one has a correspondence
 $\phi_{\ast} := p_{1\ast}\circ\psi_{\ast}:
   A_{\ast}(X^{\prime\prime})^{\circ} \vdash A_{\ast}(X)$.
(Again, there is no rational map $\phi$.)
By construction, pairs of cycle classes in
 $A_{\ast}(X^{\prime\prime})\times A_{\ast}(X)$ that are related by
 $\pi_{1\ast}\circ\psi_{\ast}$ can be extended to pairs of relative
 cycle classes in
 $A_{\ast}(W^{ss}/{\Bbb A}^1)\times A_{\ast}(W/{\Bbb A}^1)$
 with the common induced class on $\widetilde{X}$ disjoint from
 the other component in $W^{ss}_0$ or in $W_0$.

Recall the minimal lifting (curve class) $\tilde{\beta}^{\,0}$
 in $H_2(\widetilde{X};{\Bbb Z})$ of a curve class
 $\beta\in H_2(X^{\prime\prime};{\Bbb Z})$
 in a conifold degeneration from Sec.~2.

\bigskip

\noindent
{\bf Notation.} {\rm
Let $\alpha\in A_{\ast}(X^{\prime\prime})^{\circ}$.
We will adopt an abuse of notation that
 any element in the image set of the correspondence
  $\psi_{\ast}:A_{\ast}(X^{\prime\prime})^{\circ}
               \vdash A_{\ast}(\widetilde{X})$
 (resp.\ $\phi_{\ast}: A_{\ast}(X^{\prime\prime})^{\circ}
                          \vdash A_{\ast}(X)\,$)
 will be denoted by $\psi_{\ast}\alpha$
 (resp.\ $\phi_{\ast}\alpha\,$).
Since there is no chance of confusion,
 the map on the set of curve classes in
 $H_2(X^{\prime\prime};{\Bbb Z})$ to the set of curve classes in
 $H_2(X;{\Bbb Z})$ by $\beta\mapsto p_{1\ast}\tilde{\beta}^{\,0}$
 will be denoted also by $\phi_{\ast}$.
} 

\bigskip

\begin{flushleft}
{\bf GW-invariants of $X^{\prime\prime}$ in terms
     of relative GW-invariants of $(\widetilde{X},E)$.}
\end{flushleft}
Recall the moduli stack ${\mathfrak M}({\mathfrak W}^{ss},(g,n;\beta))$
 and the stabilized set $\Omega_{(g,n;\beta)}$ of admissible triples
 in Sec.~2.
Let $(\Gamma_1,\Gamma_2,I)\in \Omega_{(g,n;\beta)}$
 with $\beta_i:=b(\Gamma_i)$ and
 $\Gamma_i$ having $n_i$ legs and $r$ roots for $i=1,\,2\,$.
Then, since $K_Q=-3E$, the additivity relation of the (complex)
 dimension of virtual fundamental classes in Sec.~3.1 now reads:
 $$
  -\beta\cdot K_{X^{\prime\prime}} + n \;
 =\;   \left(\, -\beta_1\cdot K_{\widetilde{X}} + n_1
                + r - \beta_1\cdot E\, \right)\,
   +\, \left(\, 2\beta_2\cdot E + n_2 + r\, \right)\,
   -\, 2r\,,
 $$
 where
  $\beta_1\ne 0$ unless $\beta=0$,
  $n=n_1+n_2$, and
  $\beta_1\cdot E=\beta_2\cdot E\ge r\ge 0$.

There are two cases.

\bigskip

\noindent
{\it Case $(a)\,$}:
{\it $\vdim\overline{\cal M}_{g,\,0}(X^{\prime\prime},\beta)=0$.}
 In this case, one has the identity
   $$
    0\;
    =\;   \left(\, -\beta_1\cdot K_{\widetilde{X}} + r - \beta_1\cdot E\,
          \right)\,
      +\, \left(\, 2\beta_2\cdot E + r \right)\,
      -\, 2r\,,
   $$
  which implies that
  $\beta_1\cdot K_{\widetilde{X}}=\beta_1\cdot E$.
 Thus,
   $$
    \vdim {\mathfrak M}(\widetilde{Y},E;\Gamma_1)\;
     =\; r-2\beta_1\cdot E
     \hspace{1em}\mbox{and}\hspace{1em}
    \vdim {\mathfrak M}({\widetilde{Y},E;\Gamma_2})\;
     =\; 2\beta_2\cdot E + r\,.
   $$
 This implies that, if $\beta_2\ne 0$, then
  $\vdim {\mathfrak M}(\widetilde{X}, E;\Gamma_1)<0$.
 Consequently, the admissible triples in $\Omega_{(g,0;\beta)}$
  whose corresponding summand in the degeneration formula in Sec.~2
  contributes must have $\beta_2=0$ and
 the degeneration formula in Sec.~2 reads
  $$
   \Psi^{X^{\prime\prime}}_{(g,0;\,\beta)}\;
    =\; \sum_{\scriptsize
           \begin{array}{c}
             l_{1,1}+l_{1,2} = \tilde{\beta}^{\,0}\cdot E \\
             l_{1,1},\, l_{1,2}\in {\scriptsizeBbb Z}_{\ge 0}
           \end{array}
             }\,
          \Psi^{(\widetilde{X},E)}
           _{(g,0;\,\tilde{\beta}^{\,0}+l_{1,1}[C]+l_{1,2}[C^{\prime}])}\,,
  $$
  where recall that
   $\tilde{\beta}^{\,0}\in H_2(\widetilde{X};{\Bbb Z})$
    is the minimal lifting of $\beta$ in Sec.~2 and that 
   $E=C\times C^{\prime}$ from the diagram for the standard flop
    in Sec.~3.1.
 For later use, we denote the (finite) set of $l_{1,1}$-values
  in the above summation by $I_{\beta}$.

\bigskip

\noindent
{\it Case $(b)\,$}:
{\it $\vdim\overline{\cal M}_{g,0}(X^{\prime\prime},\beta)>0$.}
In this case, one considers the Gromov-Witten invariants
 $\Psi^{X^{\prime\prime}}_{(g,n;\,\beta)}(\alpha_1,\,\cdots,\,\alpha_n)$
 with the summation of the codimension of $\alpha_i$
 equal to $\vdim\overline{\cal M}_{g,n}(X^{\prime\prime},\beta)$.
Since
 $\Psi^{X^{\prime\prime}}_{(g,n;\,\beta)}
                     (\,[X^{\prime\prime}]\,,\,\cdots\,)
  =0=\Psi^X_{(\cdots)}([X],\,\cdots\,)$
 and our final goal is to relate
 $\Psi^{X^{\prime\prime}}_{(g,n;\,\beta)}
  (\alpha_1,\,\cdots,\,\alpha_n)$
 to Gromov-Witten invariants of $X$,
we will assume
  that $\alpha_i\in A_{\ast}(X^{\prime\prime})^{\circ}$ and
  that $\alpha_i\ne [X^{\prime\prime}]$.
Recall the construction of the classes
 $\psi_{\ast}\alpha_i\in A_{\ast}(\widetilde{X})$.
As cycles on $W^{ss}_0$, they are all disjoint from $Q$
 by construction.
This implies that the only admissible triples $(\Gamma_1.\Gamma_2,I)$
 in $\Omega_{(g,n;\beta)}$ that contribute to the summation in
 the degeneration formula in Sec.~2 must have $\beta_2:=b(\Gamma_2)=0$
 as well.
Thus, we have a similar expression as in Case (a):
 $$
  \Psi^{X^{\prime\prime}}_{(g,n;\,\beta)}(\alpha_1,\,\cdots,\,\alpha_n)\;
   =\; \sum_{\scriptsize
          \begin{array}{c}
            l_{1,1}+l_{1,2} = \tilde{\beta}^{\,0}\cdot E \\
            l_{1,1},\, l_{1,2}\in {\scriptsizeBbb Z}_{\ge 0}
          \end{array}
            }\,
         \Psi^{(\widetilde{X},E)}
          _{(g,n;\,\tilde{\beta}^{\,0}+l_{1,1}[C]+l_{1,2}[C^{\prime}])}
          (\psi_{\ast}\alpha_1,\,\cdots,\,\psi_{\ast}\alpha_n)\,.
 $$
Again, we denote the (finite) set of $l_{1,1}$-values
 in the above summation by $I_{\beta}$.

\bigskip

\begin{flushleft}
{\bf Gromov-Witten invariants under a small extremal transition.}
\end{flushleft}
We now want to convert the previous expression of
 $\Psi^{X^{\prime\prime}}_{(g,n;\,\beta)}(\alpha_1,\,\cdots,\,\alpha_n)$
 to one that is recognizable as a summation of Gromov-Witten
 invariants $\Psi^X_{(\cdots)}(\cdots)$ of $X$.
Recall the blow-up $p_1: \widetilde{X}\rightarrow X$.
The curve class $\tilde{\beta}^{\,0}+l_{2,1}[C]+l_{2,2}[C^{\prime}]$
 is the unique lifting on $\widetilde{X}$ of
 the curve class $\tilde{\beta}^{\,0}+l_{2,1}[C]$ on $X$
 that satisfies the condition that
  $(\tilde{\beta}^{\,0}+l_{2,1}[C]+l_{2,2}[C^{\prime}])\cdot E=0$.
Observe that
 $\tilde{\beta}^{\,0}+l_{2,1}[C]$ can never be a pure multiple of
  $[C]$ unless $\beta=0$
 and that
  $\psi_{\ast}\alpha_i
    ={p_1^{-1}}_{\ast}(\phi_{\ast}\alpha_i)$
  by construction.
It follows thus from the discussion in Sec.~3.1 that
 $$
  \Psi^{(\widetilde{X},E)}
     _{(g,n;\,\tilde{\beta}^{\,0}+l_{1,1}[C]+l_{1,2}[C^{\prime}])}
     (\psi_{\ast}\alpha_1,\,\cdots,\,\psi_{\ast}\alpha_n)\;
  =\; \Psi^X_{(g,n;\,\phi_{\ast}(\beta)+l_{1,1}[C])}
       (\phi_{\ast}\alpha_1,\,\cdots,\phi_{\ast}\alpha_n)\,.
 $$
This implies that:

\bigskip

\noindent
{\bf Theorem 3.2.1 [GW-invariant under small extremal transition].}
{\it
 With notations from above, for $\beta\ne 0$,
  $$
   \Psi^{X^{\prime\prime}}_{(g,n;\,\beta)}
     (\alpha_1,\,\cdots,\,\alpha_n)\;
    =\; \sum_{l \in I_{\beta}}\,
         \Psi^X_{(g,n;\,\phi_{\ast}(\beta)+l[C])}
           (\phi_{\ast}\alpha_1,\,\cdots,\,\phi_{\ast}\alpha_n)\,.
  $$
 If $X^{\prime\prime}$ is a Calabi-Yau $3$-fold,
  then the above identity reduces to
  $$
   \Psi^{X^{\prime\prime}}_{(g,0;\,\beta)}\;
    =\; \sum_{l \in I_{\beta}}\,
         \Psi^X_{(g,0;\,\phi_{\ast}(\beta)+l[C])}\,.
  $$
} 


\noindent
{\it Remark 3.2.2 $[\,$map on $H_2\,]$.}
 Topologically there is a surjective map (with the notation of [L-R])
  $\varphi_e: H_2(X;{\Bbb Z})\rightarrow H_2(X^{\prime\prime};{\Bbb Z})$
  from the composition of canonical homomorphisms
  $H_2(X;{\Bbb Z})
    \longrightaarrow H_2(\underline{X};{\Bbb Z})
    \stackrel{\sim}{\longrightarrow} H_2(X^{\prime\prime};{\Bbb Z})$.
 The kernel of $\varphi_e$ is generated by $[C]$.
 By construction, $(\varphi_e\circ\phi_{\ast})(\beta)=\beta$.

\bigskip

For $\beta\ne 0$, the finite set of curve classes
 $\{\,\phi_{\ast}(\beta)+l[C]\,|\, l\in I_{\beta}\,\}$
 only embeds in $\varphi_e^{-1}(\beta)\,$,
 and the latter can contain other curve classes as well.
However, for any $\hat{\beta}$ in
 $\varphi_e^{-1}(\beta)
  -\{\phi_{\ast}(\beta)+l[C]\,|\, l\in I_{\beta}\}$,
if $\hat{\beta}$ is not a curve class,
then
 $\Psi^X_{(g,n;\,\hat{\beta})}
           (\phi_{\ast}\alpha_1,\,\cdots,\,\phi_{\ast}\alpha_n)=0$
 automatically;
if $\hat{\beta}$ is a curve class, then note that
 there is no curve class $\tilde{\beta}$ in
 $H_2(\widetilde{X};{\Bbb Z})$ that satisfies both
  $p_{1\ast}(\tilde{\beta})=\hat{\beta}$ and $\tilde{\beta}\cdot E=0$
 since all such $\hat{\beta}$ are already contained in
  $\{\phi_{\ast}(\beta)+l[C]\,|\, l\in I_{\beta}\}$
  by construction.
The assumption that $\beta\ne 0$ implies also that such
 $\hat{\beta}$ cannot be a non-trivial multiple of $[C]$.
It follows then from the discussion in Sec.~3.1 that
 all such $\Psi^X_{(g,n;\,\hat{\beta})}
            (\phi_{\ast}\alpha_1,\,\cdots,\,\phi_{\ast}\alpha_n)$
 must vanish as well
since the degeneration $W/{\Bbb A}^1$ associated to blowing up $X$
 along $C$ in Sec.~3.1 relates $\Psi^X_{(g,n;\,\hat{\beta})}$ now
 to the relative Gromov-Witten invariant
 $\Psi^{(\widetilde{X},E)}_{\Gamma}$ of $(\widetilde{X},E)$
 with $\Gamma\in \emptyset$ and the latter vanishes.
Consequently, Theorem 3.2.1 above can be re-written simplier
 in terms of a superficial infinite sum:

\bigskip

\noindent
{\bf Theorem 3.2.1$^{\mathbf \prime}$
     [GW-invariant under small extremal transition].}
([L-R: Theorem B and Corollary B.1].)
{\it
 With notations from above, for $\beta\ne 0$,
  $$
   \Psi^{X^{\prime\prime}}_{(g,n;\,\beta)}
     (\alpha_1,\,\cdots,\,\alpha_n)\;
    =\; \sum_{\hat{\beta}\in \varphi_e^{-1}(\beta)}\,
         \Psi^X_{(g,n;\,\hat{\beta})}
           (\phi_{\ast}\alpha_1,\,\cdots,\,\phi_{\ast}\alpha_n)\,.
  $$
 If $X^{\prime\prime}$ is a Calabi-Yau $3$-fold,
  then the above identity reduces to
  $$
   \Psi^{X^{\prime\prime}}_{(g,0;\,\beta)}\;
    =\; \sum_{ \hat{\beta}\in \varphi_e^{-1}({\beta})}\,
         \Psi^X_{(g,0;\,\hat{\beta})}\,.
  $$
} 

\bigskip

Consider now the $3$-point function for $X^{\prime\prime}$
 $$
  \Psi^{X^{\prime\prime}}(\alpha_1,\alpha_2,\alpha_3)\;
  =\; (\alpha_1 \cdot \alpha_2 \cdot \alpha_3)_0\;\;
      + \sum_{\beta\in H_2(X^{\prime\prime};{\scriptsizeBbb Z})-\{0\}}\,
         \Psi^{X^{\prime\prime}}_{(0,3;\,\beta)}
                  (\alpha_1,\alpha_2,\alpha_3)\,q^{\beta}\,.
 $$
By construction,
 $\alpha_1 \cdot \alpha_2 \cdot \alpha_3
  = \phi_{\ast}\alpha_1
    \cdot \phi_{\ast}\alpha_2
    \cdot \phi_{\ast}\alpha_3$
since the relative cycles
 $\overline{\alpha}_i\in A_{\ast}(W^{\prime\prime}/{\Bbb A}^1)$
 used to associate
 $\alpha_i\in A_{\ast}(X^{\prime\prime})^{\circ}$
 to $\phi_{\ast}\alpha_i\in A_{\ast}(X)$
 are by definition flat over ${\Bbb A}^1$ and their lifting
 to $W^{ss}/{\Bbb A}^1$ are disjoint from the exceptional locus $E$
 of $p_1:\widetilde{X}\rightarrow X$.
Consequently, Theorem 3.2.1$^{\prime}$ and Remark 3.2.2 imply that

{\small
\begin{eqnarray*}
 \lefteqn{
   \Psi^X(\phi_{\ast}\alpha_1,
          \phi_{\ast}\alpha_2,\phi_{\ast}\alpha_3) }\\[.6ex]
  && \hspace{1em}
   =\;
   (\phi_{\ast}\alpha_1 \cdot \phi_{\ast}\alpha_2
           \cdot \phi_{\ast}\alpha_3)_0\;\;
     + \sum_{\tiny
          \begin{array}{c}
            \hat{\beta} \in H_2(X;{\tinyBbb Z})-\{0\} \\[.6ex]
            \hat{\beta}\ne m[C]
          \end{array}
             }\;\:
        \Psi^X_{(0,3;\,\hat{\beta})}
         (\phi_{\ast}\alpha_1,\phi_{\ast}\alpha_2,
          \phi_{\ast}\alpha_3)\,q^{\hat{\beta}}  \\[.6ex]
  && \hspace{2em}
       \mbox{by noting that}\;
        \sum_{\mbox{\tiny $\hat{\beta}=m[C],\,m>0$}}
         \Psi^X_{(0,3;\,\hat{\beta}) }
         (\phi_{\ast}\alpha_1,
          \phi_{\ast}\alpha_2, \phi_{\ast}\alpha_3)\,
          q^{\hat{\beta}}\;  =\; 0\;\;
        \mbox{since $\;C\cap \phi_{\ast}\alpha_i=\emptyset$}\,,\\[.6ex]
  && \hspace{1em}
   =\;
   (\phi_{\ast}\alpha_1 \cdot \phi_{\ast}\alpha_2
          \cdot \phi_{\ast}\alpha_3)_0\;\; 
     + \sum_{\mbox{\tiny
             $\beta
              \in H_2(X^{\prime\prime};{\tinyBbb Z})-\{0\}$} }\;\:
       \sum_{\mbox{\tiny $\hat{\beta}\in \varphi_e^{-1}(\beta)$}} 
        \Psi^X_{(0,3;\,\hat{\beta})}
         (\phi_{\ast}\alpha_1,\phi_{\ast}\alpha_2,
          \phi_{\ast}\alpha_3)\,q^{\hat{\beta}} \\[.6ex]
  && \stackrel{\rule[-.6ex]{0ex}{1ex}
               q^{\hat{\beta}}\,
                \rightarrow\, q^{\varphi_e(\hat{\beta})}\,
                =\, q^{\beta}}
              {\longrightarrow}\;
   (\alpha_1 \cdot \alpha_2 \cdot \alpha_3)_0\;\;
     + \sum_{\mbox{\tiny
               $\beta
                \in H_2(X^{\prime\prime};{\tinyBbb Z})-\{0\}$}}\;\:
        \Psi^{X^{\prime\prime}}_{(0,3;\,\beta)}
         (\alpha_1,\alpha_2,\alpha_3)\,q^{\beta} \\[.6ex]
  && \hspace{1em}
   =\; \Psi^{X^{\prime\prime}}(\alpha_1,\alpha_2,\alpha_3)\,.
\end{eqnarray*}
} 

In summary:

\bigskip

\noindent
{\bf Corollary 3.2.3 [${\mathbf 3}$-point function under small extremal
     transition].} ([L-R: Corollary B.2].)
{\it
 The $3$-point function
  $\Psi^{X^{\prime\prime}}(\alpha_1,\alpha_2,\alpha_3)$ of
  $X^{\prime\prime}$
  is equal to the $3$-point function
  $\Psi^X(\phi_{\ast}\alpha_1,
          \phi_{\ast}\alpha_2,\phi_{\ast}\alpha_3)$
  of $X$ with a change of variables
   $q^{\beta}\rightarrow q^{\varphi_e(\beta)}$
   over $\beta\in H_2(X;{\Bbb Z})\,$.
} 

\bigskip

\noindent
{\it Remark 3.2.4
$[\,$K\"{a}hler $3$-fold hierarchy vs.\ GW-hierarchy$\,]$.}
The discussions of this section can be summarized schematically
by the following diagram:
 $$
  \parbox{12ex}{\raggedright\scriptsize\it
     reduction of \newline K\"{a}hler  \newline deformations}\!\!
  \mbox{\LARGE $\downarrow$}\hspace{2em}
  \begin{array}{rcc}
   X \hspace{1.4ex}
    & \stackrel{\rule[-.6ex]{0ex}{1ex}\mbox{\scriptsize\it
                flop $\, \phi$}}{\dashrightarrow}
    & X^{\prime} \\
   \parbox{10ex}{\scriptsize\it
                 small \newline extremal \newline
                 transition}\!\! \top \hspace{1.4ex}
    && \\
   X^{\prime\prime}
  \end{array}\;
  \Longrightarrow\;
  \begin{array}{ccc}
   \GW(X)  & \stackrel{\sim}{\longrightarrow} & \GW(X^{\prime}) \\
    \downarrow && \\
   \GW(X^{\prime\prime}) &&
  \end{array}.
 $$

\vspace{9em}
\begin{flushleft}
{\Large\bf Appendix. The stringy and the symplectic aspect.}
\end{flushleft}
A very terse account of transformations of Gromov-Witten
 invariants from the stringy and the symplectic viewpoint
 is given in this appendix to complement the main text
 and to provide a slightly more (though never) complete picture
 of where the issue sits in all.

\bigskip

\begin{flushleft}
{\large\bf A.1\hspace{1ex}
  Transformation of GW-invariants from the stringy viewpoint.}
\end{flushleft}
The generating function of Gromov-Witten invariants gives
 a mathematical expression/ formulation for some correlation
 functions in superstring theory.
The transformation of these correlation functions from one phase to
 another is an important issue in the theory.
In this subsection, we discuss the concept of
 {\it phase structures} in superstring theory, focusing on those
 sharing an interface with birational geometry, and conclude with
 an explanation of {\it Witten's wall-crossing formula} in [Wit5]
 and its completion in [L-R].

\bigskip

\begin{flushleft}
{\bf A glimpse of superstring theory.}
\end{flushleft}
String theory is a constantly fast-growing subject, containing
 a vast amount of contents with an ultimate goal to decipher
 God's code of creation and to unravel the mystery of
 the (physical) Universe.
So far no-one can predict if the theory will finally settle to
 a final form, yet mathematicians have already witnessed its
 horrifying power of synthesizing drastically different branches
 of mathematics into an integrated entity via stringy dualities.
Such ability and strength to bring together the microscopic world
 (particle physics), the macroscopic world (cosmology), and
 the man-made world (mathematics) make string theory very unique.
A very conservative glimpse of it is given in {\sc Diagram A-1-1},
 whose details are referred to [G-S-W], [L\"{u}-T], [Polc1],
 and the literatures guided therein up to the end of year 1997.
Surveys on more recent developments can be found in TASI Lectures
 from year 1996 on.
[Zw] contains a lot of physicists' sense/intuitions/insights
 cleanly explained, which can be very helpful for mathematicians.

\begin{figure}[htbp]
 \setcaption{\small
 {\sc Diagram} A-1-1.
  \baselineskip 12pt
  The {\it string world-sheet}, the {\it brane-world-volume},
   the {\it total space-time}, and
   the {\it $4$-dimensional effective field theory}
   aspect of a string theory.
  Scattering amplitudes of fields and D-branes in space-time
   computed via string world-sheet methods and via space-time
   field theory method have to match order by order.
  Each of the four aspects itself has a {\it Wilson's theory-space}
   associated to it, containing all the phases of the field theory
   associated to that aspect of the string theory.
  {\it Dualities} can be realized either as a local isomorphisms
   or a coordinate change on the Wilson's theory space
    ${\cal S}_{\Wilson}$
   that induces an isomorphism on the universal family of
    Hilbert spaces of states, ring of operators, and
    correlation functions of these operators over ${\cal S}_{\Wilson}$.
  In the weak string coupling regime, strings are light and branes
   are heavy and hence string are regarded as more fundamental.
  In the strong string coupling regime, branes can become light and
   strings become heavy and should be no longer treated as
   the most/unique fundamental object in the theory.
 } 
\centerline{\psfig{figure=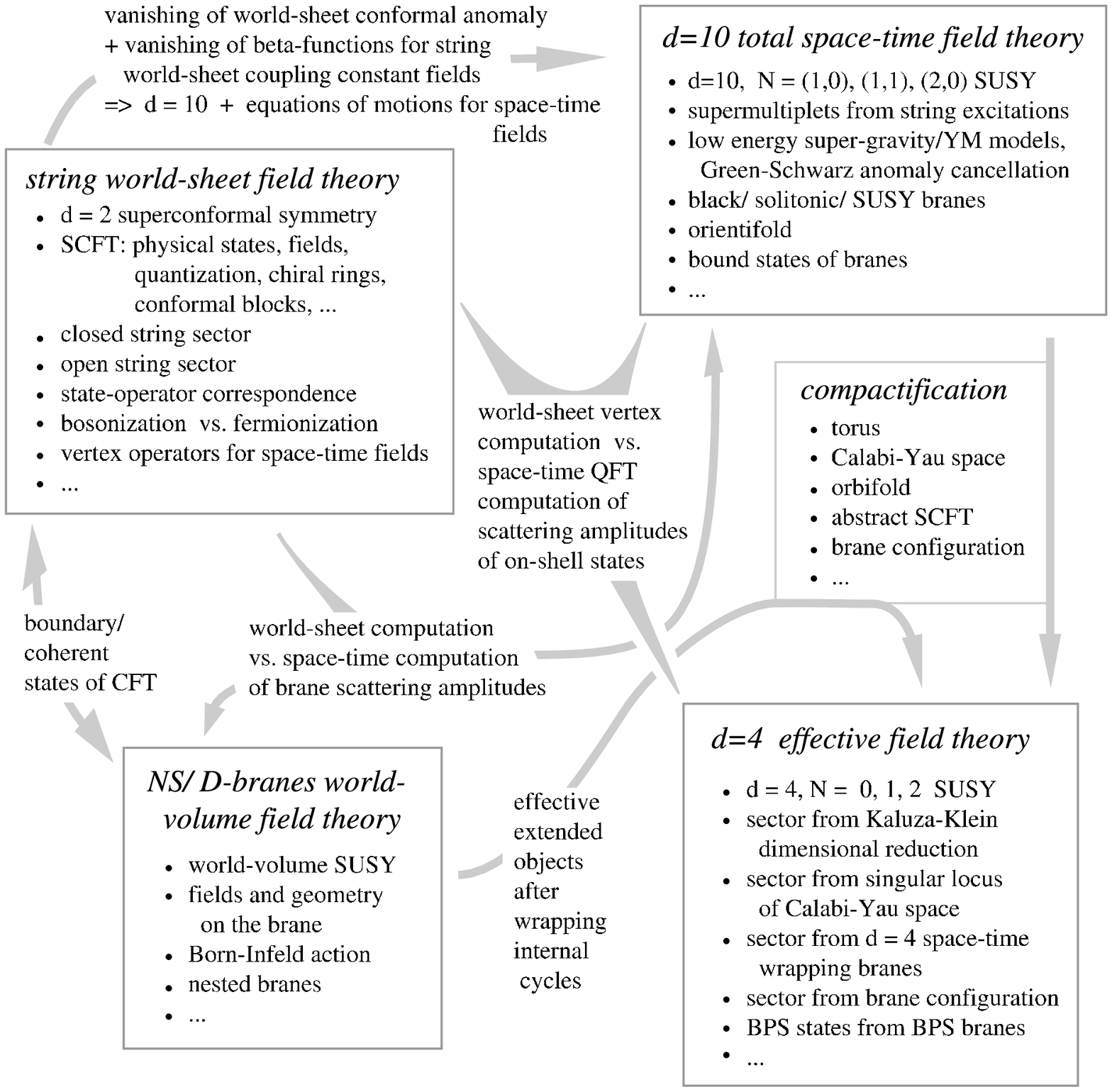,width=13cm,caption=}}
\end{figure}

\bigskip

\noindent
{\it Remark A.1.1.}
Two brief remarks follow:
\begin{itemize}
 \item[(1)]
  {\it String/M/F-theory.} 
  The five fundamental {\it superstring models} at $10$-dimensions
   are now known to be difference phases of a single {\it M-theory},
   which also contains the $11$-dimensional supergravity theory
   as one of its phases, [Wit6].
  There is also a $12$-dimensional {\it F-theory} from geometrizing
   a complexified coupling constant in the low energy supergravity
   theory associated to the $10$-dimensional IIB superstring theory,
   [Va2].
  See [Polc1: Chapter 14] for a review of various dualities that
   connect all these theories/models.
  Compactifications of these theories bring $G_2$ $7$-manifolds and
   $\Spin(7)$ $8$-manifolds into play as well.
  Transformations of correlations in superstring theory may be reduced
   to a reduction of transformations of correlations in these higher
   level theories.

 \item[(2)]
  Though {\it D-branes} had already entered string theory in late
   1980's ([D-L-P]), their very fundamental role in string
   theory began to be understood only in the mid 1990's
   ([Polc3]; see also [H-T] and [Wit6]).
  Since then, string theory has been revolutionized to such
   an extent that one has to re-ask
   {\it What is string theory?}
   {\it What is space-time?} and even 
   {\it What is a D-brane itself?}
  {\it Non-commutative geometry} is expected to play a role here,
  cf.\ [Do2].
\end{itemize}

\bigskip

\begin{flushleft}
{\bf Wilson's theory-space, phases, and dualities/correspondences.}
\end{flushleft}
A quantum field theory (QFT), as specified by its Lagrangian density,
 admits continuous deformations by, e.g.\ varying the coupling
 constants in the Lagrangian  density. {\it Wilson's theory-space},
 denoted by a general notation ${\cal S}_{\Wilson}$ in this subsection,
 is meant to be the universal parameter/moduli space that encodes all
 such deformations of a field theory with a specified combinatorial
 type like number/type of fields involved
 (e.g.\ scalor theory, gauge theory with matters, gravity theory,
        nonlinear sigma model),
  dimension and symmetry of the field theory
 (e.g.\ $d=4$, $N=2$ super-Yang-Mills theory with gauge group $\SU(n)$),
 types of couplings (e.g.\ $\phi^4$-theory, Yukawa theory), ..., etc.
This is in the same spirit as the use of parameter/moduli spaces
 in mathematics for parameterizing
 e.g.\ varieties, K\"{a}hler manifolds, coherent sheavs, bundles,
       algebras, rings, ... etc.\ with a specification of
       combinatorial type like Hilbert polynomials,
       the number of generators, ..., etc..
With these said, however, quantum field theory, quantities defined
 therein and techniques used to draw conclusions remain full with
 surprises and mysteries for mathematicians
 (e.g.\ [AH-C-G] as a more recent example).
While trying to spell out ${\cal S}_{\Wilson}$ below, we have to
 admit as well our ignorance to ${\cal S}_{\Wilson}$ and the universal
 family, in notation ${\cal U}_{\Wilson}$,
 that goes with it (as in any moduli problem in mathematics).
See e.g.\ [Ca], [P-S], [Poly], [Stra], [A-D-S], [G-V-W], [Polc2]
 for physicists' insights and [B-B] for a discussion on the
  mathematical structure in renormalization of QFT.
\begin{itemize} 
 \item[$\bullet$]
 {\it ${\cal S}_{\Wilson}$ as a universal parameter space.}
 As a universal parameter space for QFT of a specified combinatorial
  type, a local coordinate system on the Wilson's theory-space
  ${\cal S}_{\Wilson}$ should contain
   \begin{itemize}
    \item[(1)]
    {\it coupling constants} of the field theories in the specified
    combinatorial types,

    \item[(2)]
    a {\it length/energy scale $\lambda$ } that sets the effectiveness
    of the theory so that any effect from a smaller smaller scale
    (or larger energy) than the specified one is regarded as already
    integrated out,

    \item[(3)]
    other {\it natural parameters} in the theory depending the contents,
    e.g.\ string tension $\alpha^{\prime}$ or
    the Planck's constant $\hbar$,

    \item[(4)]
    all the additional {\it parameters} arising {\it from
    renormalization/regularization schemes} that are needed to
    render physical quantities of the theory finite.
  \end{itemize}
 {\it Renormalization group flow} (RG-flow) on ${\cal S}_{\Wilson}$
 takes a theory in one length scale $\lambda$ to another theory in
 a larger length scale $\lambda^{\prime}>\lambda$ by integrating out
 further the effect between $\lambda$ and $\lambda^{\prime}$.

 \item[$\bullet$]
 {\it The universal object $\,{\cal U}_{\Wilson}$
      over ${\cal S}_{\Wilson}$.}
 A quantum field theory contains a long list of contents, e.g.\ 
  the domain space on which fields are defined,
  the target spaces fields take value,
  bundles/sheaves whose sections give the fields,
  the Lagrangian density,
  symmetries,
 Hilbert space of states,
  observables/operators, 
  correlation functions of fields/observables/operators,
    vacuum space/manifold $\cdots\,$.
 These pile together and form several types of universal objects
  over ${\cal S}_{\Wilson}$.
 We will denote these universal objects collectively
  by ${\cal U}_{\Wilson}$ over ${\cal S}_{\Wilson}$.

 \item[$\bullet$]
 {\it The topology and geometry on ${\cal S}_{\Wilson}$.}
  The use of coupling constants in Lagrangian densities and
   other parameters as the local coordinates on ${\cal S}_{\Wilson}$
   gives a topology on ${\cal S}_{\Wilson}$.
  Infinitesimal deformations of a QFT $T$ can be realized as
   operators $O_i$, $i\in$ some index set $I$, in the theory $T$.
  The $2$-point functions $\langle O_iO_j\rangle_T$ induces
   the (Zamolodchikov) metric on ${\cal S}_{\Wilson}$.
  This makes ${\cal S}_{\Wilson}$ a geometric object.

 \item[$\bullet$]
 {\it Phase structure on ${\cal S}_{\Wilson}$.}
 In general, ${\cal S}_{\Wilson}$ has a stratification determined
  by the behavior of the field theories it parameterizes.
 E.g.\ the $2$-point (correlation) function can behavior
  asymptotically following a power law on on one region of
  ${\cal S}_{\Wilson}$ while asymptotically exponentially
  on another region;
  or the vacuum space can be described by one deformation family of
   geometries/manifolds on one region of ${\cal S}_{\Wilson}$
  while by a different deformation family of geometries/manifolds
   on another region od ${\cal S}_{\Wilson}$.
 This gives a {\it phase structure} on ${\cal S}_{\Wilson}$
  with each maximal stratum of the stratification called
   a {\it phase} of ${\cal S}_{\Wilson}$
  (or of the theory of the specified combinatorial type).
 Some nature of the quantum field theries parameterized by
  ${\cal S}_{\Wilson}$ changes when passing from one phase to
   another.
\end{itemize}

Many mathematical moduli spaces/stacks with their universal
 object appear as a sub-fibration of ${\cal U}_{\Wilson}$
  over ${\cal S}_{\Wilson}$.
For example,
 in nonlinear sigm models
 (resp.\ compactifications of a superstring model)
 the moduli space of the target space (resp.~ compactifying space)
 (possibly with decorations like bundles/coherent-sheaves thereover
  or cycles therein)
 is naturally contained in the related ${\cal S}_{\Wilson}$.
And for gauge theories, the natural coordinates on ${\cal W}_{\Wilson}$
 can contain gauge instanton corrections, which involve integrations
  over the moduli space of connections or stable sheaves.
Such links of ${\cal S}_{\Wilson}$ to mathematical moduli problems
 turn many physical/stringy computations/statements to surprising
 mathematical conjectures,
(cf.\ the discussion on geometric engineerings below).

\bigskip

\noindent
{\it Remark A.1.2 $[$reduced theory-space$]$.}
In the physicists' setting in describing Wilson's theory-space,
 a local chart thereon consists only of the coupling constants
 (the first in our four sets of parameters).
Though physicists have developed techniques to use it fluently,
 mathematically this reduced description is less satisfying.
The distinction is similar to (but more involved than)
 the distinction between a coarse moduli space versus
 the moduli stack in a moduli problem.
However, in all our discussions below it is really this part
 that matters.
Thus, we will resume physicists' conventional use/definition
 below; i.e.\ ${\cal S}_{\Wilson}$ will be the parameter space
 of Largrangian densities of a fixed combinatorial type
 for the rest of the discussions.
Depending on the level of the theory one explores
 - classical, quantum, or non-perturbative -,
locally on ${\cal S}_{\Wilson}$ these coupling coefficients
 as coordinates may be subject to corrections at each level.
The length scale $\lambda$ at which the theory becmes effective
 comes in usually via an energy cutoff in the
 renormalization/regularization procedure.
(It can come into play via other mechanism,
  e.g.\ dynamical generation, as well.)

\bigskip

\noindent
{\bf Example A.1.3 [natural coordinates on ${\cal S}_{\Wilson}$
     and stringy duality].}
Compactifications of different superstring models on different
 Calabi-Yau $3$-folds $X$ (possibly decorated with supersymmetric
 cycles) may give rise to isomorphic $4$-dimensional effective
 supersymmetric field theories.
(Such effective field theory will be denoted in Candelas-type notation
 by $\String[X]$ in this subsection.)
When this happens, each choice of the superstring models and
 the topological type of comapctifying spaces $X$ gives rise to
 a natural local coordinate system on the Wilson's theory-space
 ${\cal S}_{\Wilson}$ of the related $d=4$ field theory.
The local coordinate transformation between these different coordinate
 systems on ${\cal S}_{\Wilson}$ is the much-explored ``mirror map".
Unfamiliar readers can use keyword search:
 ``Calabi-Yau", ``mirror symmetry", ``special geometry",
 and ``$N=1$ or $N=2$ supersymmetry" to find more detailed explanations;
 see also [C-dlO-G-P], [Gre], [Polc1: Chapter 17 and Chapter 18].
We will give another example below when discussing a link between
 birational geometry and geometric engineering of quantum field theory.

\bigskip

\begin{flushleft}
{\bf Phases in superstring theory and birational geometry.}
\end{flushleft}
Three themes in superstring theory that involve some key ingredients
 of Mori's program in birational geometry ([Kol-M]) are listed below.
(These had come before the language and techniques of derived
  categories, $t$-structures, and perverse sheaves came to play a role
  as well on both sides;
 see [As] and [Wa] for a review of these later related directions.)

\begin{itemize} 
\item[$\bullet$]
{\it Phases of string world-sheet theory vs.\
     transitions of Calabi-Yau manifolds.}
This is related to the theme of this note and we will come back to it.
See also [Stro2] for a review from the aspect of $d=4$ supersymmetric
 effective field theories via compactifications of the superstring
 models on Calabi-Yau spaces.
\end{itemize} 

\begin{itemize} 
\item[$\bullet$]
{\it Special loci in the theory-space of $d=4$ effective SQFT
     vs.\ singularities of Calabi-Yau spaces}$\,:$
{\it geometric engineering} ([K-K-V]).
Superstring theory contains also other dynamical extended objects,
 particularly D-branes.
When a superstring model is compactified on a singular Calabi-Yau
 $3$-space $X_0$ that arises from deforming either the complex or
 the K\"{a}hler structure of smooth Calabi-Yau manifolds $X$,
D-branes that are wrapped around vanishing cycles or exceptional
 locus in $X$ associated to the complex or the K\"{a}hler degeneration
 (which contribute to massive spectrum of the $d=4$ effective theory
  $\String[X]$) can give rise to extra degeneracy of the massless
  spectrum of the low-energy $d=4$ effective supersymmetric field
  theory $\String[X_0]$, cf.~[Stro1].
In particular, there can be enhanced gauge symmetry
 determined by the singular locus $X_{0\,\sing}\subset X_0$
 when one moves in the theory-space of the $d=4$ SQFT from
 $\String[X]$ to $\String[X_0]$, e.g.~[K-M-P].
Such phenomena are indeed required in establishing,
 e.g.\ IIA/Heterotic string duality and they go up also
 to compactifications of  M- and F-theory.
Focusing only on a neighborhood of $X_{0\,\sing}\subset X_0$
 and its resolution gives rise to a corrrespondence
 (geometric engineering of SQFT):
$$
 \begin{array}{ccc}
 & \parbox{14ex}{\tiny
    compactification of \newline
    a superstring model  \newline on $X_0$}        \\[-3ex]
  \framebox[27ex]{\hspace{1ex}\parbox{26ex}{\footnotesize
     local singular \newline Calabi-Yau $3$-space $X_0$ \newline
     with degeneration $X\leadsto X_0$   }}
   & \Longrightarrow
   & \framebox[31ex]{\hspace{1ex}\parbox{30ex}{\footnotesize
       $\String[X_0]\,$:  \newline
       $d=4$ SUSY gauge field theory \newline
       with the gauge group and \newline
       the supermultiplets \newline
       determined by $X_{0\,\sing}\subset X_0$ }}  \\
   &&
 \end{array}
$$

\item[]\hspace{2.8ex}
In particular, this implies that there are {\it three} natural but
 {\it different coordinates} on a {\it same region} in the Wilosn's
 theory-space ${\cal S}_{\Wilson}$ for the related $d=4$
 supersymmetric gauge theory:
 \begin{itemize}
  \item[(1)]
   from $\String[X_0]\,$: system of coordinates that contain
    string world-sheet instantons corrections;

  \item[(2)]
   from local mirror to $X\,$: system of coordinates that arise
    from periods of a (complex) curve $C$, the local mirror to $X$;

  \item[(3)]
   as a $d=4$ super Yang-Mills theory in its own right:
    system of coordinates that contain gauge instanton correcions. 
 \end{itemize}
The related local coordinate transformations on ${\cal S}_{\Wilson}$,
 phrased in mathematical language, say that some limit of generating
 function of Gromov-Witten invariants for the local Calabi-Yau
 $3$-fold $X$ can be related to the generating function of
 Donaldson/Seiberg-Witten invariants for the $d=4$ gauge theory
 (on an appropriate compactification of ${\Bbb R}^4$ or ${\Bbb C}^2$)
 and that both can be reproduced by periods of curves;
see [N-Y] for some related review.
On the other hand, a direct relation at the level of the moduli spaces
 involved (with the virtual fundamental classes thereon constructed)
 remains very technical to understand.

\vspace{-1.6ex}
\item[]\hspace{2.8ex}
Incidentally, a singular $X_0$ can have different resolutions
 $X\leadsto X_0$.
Thus these spacial loci can be where different phases/branches
 of the theory-space meet.
\end{itemize} 

\begin{itemize} 
\item[$\bullet$]
{\it Phases of D-brane world-volume theory vs.\
     resolutions of quotient singularities.}
A D-brane carry charges for RR-fields in type II string theories.
It is also where end-points of open strings are sticked to.
When a D$3$-brane sits in the singular locus of
 a quotient space-time ${\Bbb R}^{3+1}\times {\Bbb C}^3/\Gamma$,
 where $\Gamma$ is a finite subgroup in $U(3)$ and acts on
 ${\Bbb C}^3$ effectively and freely except at the origin,
the induced field theory on the world-volume of the D-brane
 is a $d=4$ SQFT with supermultiplets from dimensional reduction,
 orbifolding, and also excitations of open strings with one or
 both end-points on the D-brane.
The contents of this $d=4$ field theory can be encoded in a quiver
 diagram.
Similar to the discussion in [Wit5], there are coupling constants
 in the theory that parameterizes the theory-space
 ${\cal S}_{\Wilson}^{\,\mbox{\it\tiny D-brane}}$
 of such field theories on the D-brane.
On different regions of ${\cal S}_{\Wilson}^{\,\mbox{\it\tiny D-brane}}$
 the associated vacuum manifold of the theory can be different,
 giving rise to a phase structure on
 ${\cal S}_{\Wilson}^{\,\mbox{\it\tiny D-brane}}$.

\vspace{-1.6ex}
\item[]\hspace{2.8ex}
As a probe to the nature of its ambient space-time,
before descended to a field theory expanded at a vacuum manifold,
the field theory on the D-brane,
 regarded as a $d=4$ nonlinear sigma model with target
  the transverse space to the brane,
 ``sees" a non-commutative ambient space with local coordinates
 transverse to the brane enhanced to Lie-algebra-valued depending
 on $\Gamma$.
However, the vacuum manifolds obtained from solving minimal potential
 equation with flatness conditions from the superpotential
 recover the commutative nature of the transverse geometry to the brane
 yet they do not go back to the singular ${\Bbb C}^3/{\Gamma}$
 that one starts with.
Rather, these vacuum manifolds are various (partial) resolutions of
 ${\Bbb C}^3/\Gamma$.
Thus, moving from one phase to another in
 ${\cal S}_{\Wilson}^{\,\mbox{\it\tiny D-brane}}$ and
 via the descendant low energy field theory on the D-brane
 world-volume, the D-brane ``sees" different ordinary
 (i.e.\ commutative) ambient space-times from the resolutions of
 ${\Bbb C}^3/\Gamma$, e.g.\ [D-M], [D-G-M], and [G-L-R].
In this way, the phase structure on the induced field theory
 on the D-brane world-volume is linked to resolutions of
 a quotient singularity. {\sc Figure} A-1-2.

\vspace{-1.6ex}
\item[]\hspace{2.8ex}
(It should be noted that
 there are other perspectives of phases of $D$-branes;
  see e.g.\ [H-I-V] and [Ma] in connection with
  McKay correspondence, mutations, and helixes; and
  e.g.\ [A-H], [H-W], and [Wit8] for phase structures of
  field theory from brane configurations.)
\end{itemize} 

\begin{figure}[htbp]
 \setcaption{\small
 {\sc Figure} A-1-2.
  \baselineskip 12pt
   The D-brane sitting along the singular locus perceives 
    a noncommutative thickening of the singular space-time
    to begin with.
   At low energy, the space-time it perceives resumes 
    an ordinary yet different space-time: a resolution of 
    the original singular space-time.
   Each resolution corresponds to a geometric phase in
    the theory-space
    ${\cal S}_{\Wilson}^{\,\mbox{\it\tiny D-brane}}$
    of the QFT on the D-brane world-volume.
 } 
\centerline{\psfig{figure=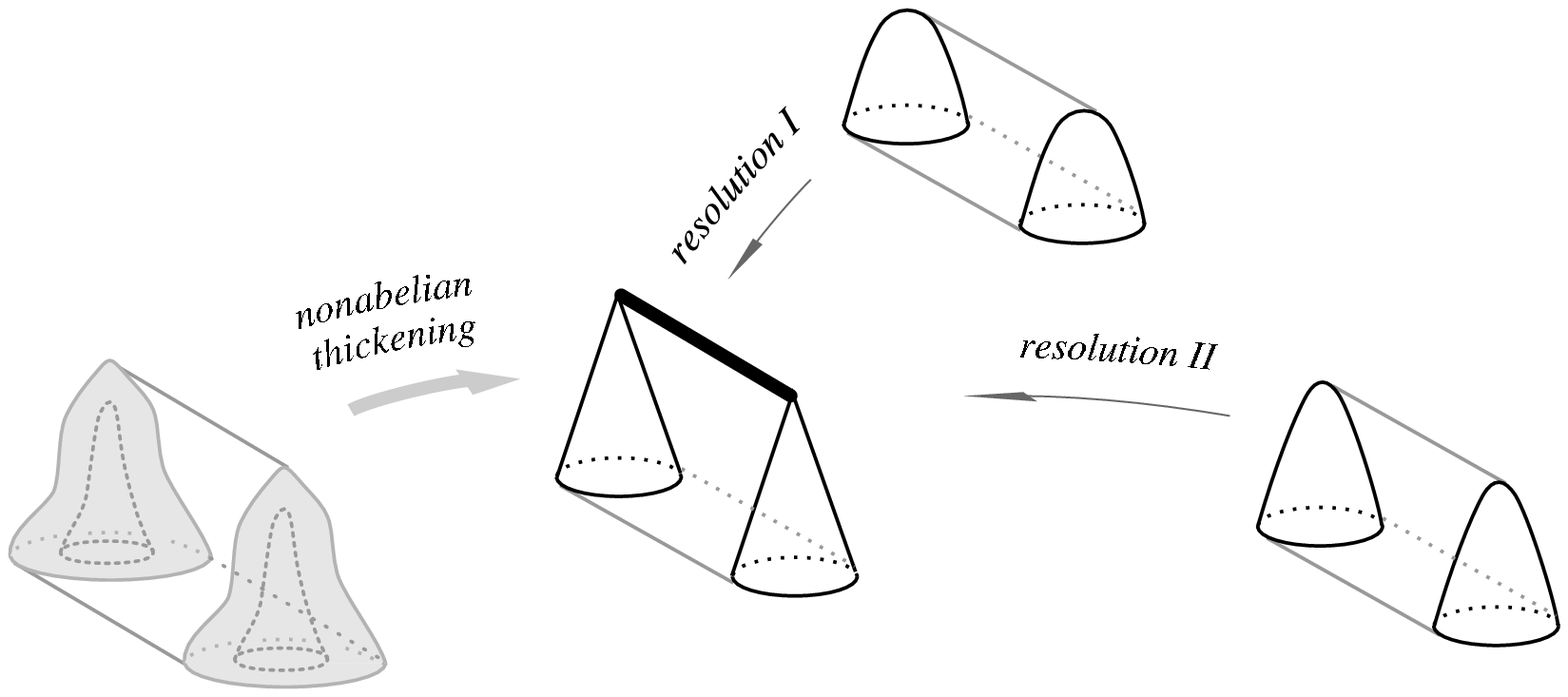,width=13cm,caption=}}
\end{figure}

The world-sheet aspect of a superstring theory leads directly
 to the Gromov-Witten theory. We will now focus on this link.
Some related physics background is highlighted in
 [Liu-L-Y: Appendix], which we will use directly.

\bigskip

\begin{flushleft}
{\bf From superstrings to Gromov-Witten theory: A-model.}
\end{flushleft}
\nopagebreak[4]
The path-integral quantization of superstrings in principle has
 the moduli spaces of holomorphic maps hidden in its integral.
However, the technical issue of existence of global spinors on
 general string world-sheets and the fact that the path-integral
 involves more than such maps make these moduli spaces at best hidden.
In [Wit4], Witten introduces the notion of topological twists.
Associated to a Calabi-Yau manifold $X$ are two $d=2$, $N=(2,2)$
 topological field theories: the A-model and the B-model for $X$.
And the path-integral of the A-model $\Amodel[X]$ for $X$ localizes
 on the space of holomorphic maps into $X$.
This brings the moduli space of holomorphic maps of arbitrary
 genus manifestly out in superstring theory
and links superstring theory directly to Gromov-Witten theory.

\bigskip

\begin{flushleft}
{\bf Phase structure and
     the transformation of A-model correlation functions.}
\end{flushleft}
Though one's final goal is to understand the $d=4$ effective
 theory $\String[X]$ from superstring compactification on
 a Calabi-Yau $3$-space $X$,
some quantities in $\String[X]$ can be computed via
 the $d=2$ world-sheet theory $\Amodel[X]$. 
In particular, one has the following diagram of
 local embeddings
$$
 \begin{array}{cccl}
  {\cal S}_{\Wilson}^{\,d=2,\,N=(2,2)}
    & \stackrel{\rule[-.6ex]{0ex}{1ex}\mbox{\scriptsize\it
                A-model$\,[\,\cdot\,]$}}{\longleftarrow}
   & \framebox[24ex]{\parbox{21ex}{\it\small
       Moduli space ${\cal M}_{CY^3}$ of \newline
       $\hspace{1ex}$ Calabi-Yau $3$-spaces }}  \\[2.6ex]
  && \hspace{7ex}\downarrow
     \mbox{\scriptsize\it String$\,[\,\cdot\,]$}\\[1ex]
  && {\cal S}_{\Wilson}^{\,d=4,\,N=2} &,
 \end{array}
$$
where the horizontal arrow is the assignment to $X$ the
 $2$-dimensional field theory $\Amodel[X]$ and the vertical arrow
 is the assignment to $X$ the $4$-dimensional field theory $\IIA[X]$
 in Candelas' notation (one may use $\IIB[X]$ as well).
(The image of the horizontal arrow indeed lies in a small subspace
  of ${\cal S}_{\Wilson}^{\,d=2,\,N=(2,2)}$ that parameterizes
  theories with also conformal symmetries.)
Moving from one irreducible component to another irreducible component
 of ${\cal M}_{CY^3}$ by deforming the complex and/or the K\"{a}hler
 structure of the $3$-space corresponds to moving from one phase
 to another phase in both
 ${\cal S}_{\Wilson}^{\,d=2,\,N=(2,2)}$  and
 ${\cal S}_{\Wilson}^{\,d=4,\,N=2}$.
How correlation functions of the field theory transform when
 crossing the wall/boundary between different geometric phases
 are natural question to ask.
For the phase transition corresponding to a standard flop of
  Calabi-Yau $3$-fold, Witten conjectured a wall-crossing
  formula in [Wit5: Sec.~5.5].
(Some notations follow [Wit5] and are used only in this subsection.
   In particular, $E$ below is not our exceptional divisor $E$
   in Sec.~3.1.
  See also the related discussions in [K-M-P].)

Consider a $d=2$, $N=(2,2)$ gauged linear sigma model (GLSM)
 with a fixed number of vector multiplets, chiral multiplets
 and appropriate choice of charge vectors so that
the vacuum manifold in each geometric phase of the theory is
 a Calabi-Yau $3$-fold.
The descendant $d=2$ non-linear sigma models associated to
 the geometric phases of the GLSM lies in the image
 $\Amodel[{\cal M}_{CY^3}]$;
they can also be identified as a subspace of $\String[{\cal M}_{CY^3}]$
 in ${\cal S}_{\Wilson}^{\,d=4,\,N=2}$.
 \begin{itemize}
  \item[(1)]
   Following [Wit5: Sec.~5.5], suppose that $X_+ := X$ and
    $X_- := X^{\prime}$ are realized as the respective vacuum
    manifold in two adjacent geometric phases of the GLSM.
   On each such phase, there is a descendant non-linear sigma model
    with target the vacuum manifold $X^{\pm}$.
   The (real) $2$-dimensional space on which fields of
    the non-linear sigma models are defined can be identified
    as the (Wick-rotated) string world-sheet.
   Let $\alpha_1=\alpha_2=\alpha_3=$ a divisor $E$ on $X_+$ that
    intersects $C$ properly with $[C]\cdot E=+1$.
   They correspond to observables $O_E$ of the $d=2$ nonlinear
    sigma model with target $X_+$.
   Then under the topology change $X_+\dashrightarrow X_-$
    (which is our map $\phi$), $[C^{\prime}]\cdot E=-1$.
   Let $\lambda_+$ (resp.\ $\lambda_-$) be the contribution to
    the $3$-point function of the non-linear sigma model
    $\langle O_E\, O_E\, O_E\rangle_{X_+} =\Psi^{X_+}(E,E,E)$
    (resp.\ $\langle O_E\, O_E\, O_E\rangle_{X_-}=\Psi^{X_-}(E,E,E)$,
     which means $\Psi^{X_-}(\phi_{\ast}E,\phi_{\ast}E, \phi_{\ast}E)$
     in Sec.~3.1)
    purely from the world-sheet instantons wrapped around $C$
     (resp.\ $C^{\prime}$).
   Then, following [C-dlO-G-P] and [A-M], Witten deduced that
    $$
     \mbox{$\hspace{12em}$}
     \lambda_+ - \lambda_- \; =\; -1\,.
     \hspace{6em}\mbox{[Wit5: Eq.(5.48)]}
    $$
   The ``=" here is a direct formal manipulation on power series
    and can be interpreted as an analytic continuation from
    $\lambda_++1$ to $\lambda_-$.
   Witten argued that it cannot be just $\lambda_+=\lambda_-$ as
    one might naively expect and that the discrepancy $-1$ in
    the above equation has to do with the classical $3$-point
    function, namely the triple interseection product
     ($=E\cdot E\cdot E$ in this example) of the observables
     (i.e.\ divisors on $X_{\pm}$) chosen..
   Taking together the classical $3$-point function plus
    the world-sheet instanton contribution,
   Witten concludes that the new stringy corrected $3$-point function
    of the non-linear sigma model in the two phases of the gauged
    linear sigma model should be related by an analytic continuation
    when crossing the wall between the two phases.
 \end{itemize}
 \begin{itemize}
  \item[(2)]
   Compared with [L-R] and Sec.~3.1,
    the $\lambda^{\pm}$ is the third term of the decomposition
    of the full $3$-point functions
     $\Psi^X(\alpha_1,\alpha_2,\alpha_3)$ and
     $\Psi^{X^{\prime}}(\phi_{\ast}\alpha_1,
           \phi_{\ast}\alpha_2, \phi_{\ast}\alpha_3)$.
   Witten's statement above is the identity
   \begin{quote}
    {\footnotesize
     \begin{eqnarray*}
       \lefteqn{
        \phi_{\ast}\left(\,\alpha_1\cdot\alpha_2\cdot\alpha_3
         + ([C]\cdot\alpha_1)([C]\cdot\alpha_2)([C]\cdot\alpha_3)\,
           \frac{q^{[C]}}{1-q^{[C]}}\,
                     \right)            }\\
       && \stackrel{a.c.}{=}\;
               \phi_{\ast}(\alpha_1)\cdot\phi_{\ast}(\alpha_2)
                                    \cdot\phi_{\ast}(\alpha_3)
              + ([C^{\prime}]\cdot\phi_{\ast}(\alpha_1))
                ([C^{\prime}]\cdot\phi_{\ast}(\alpha_2))
                ([C^{\prime}]\cdot\phi_{\ast}(\alpha_3))\,
                \frac{q^{[C^{\prime}]}}{1-q^{[C^{\prime}]}}
      \end{eqnarray*}
      } 
   \end{quote}
   (with the necessary notation modification in [L-R] to adjust
    to the push-pull direction of differential forms).
 \end{itemize}
 \begin{itemize}
  \item[(3)]
   The remaining contributation to the full $3$-point function
    is the second term in the decomposition.
   From [L-R] and Sec.~3.1, it is intact;
   namely, it follows simply by applying $\phi_{\ast}$ without
    having to perform analytic continuations.
   In this aspect, the only difference for $3$-point functions
    when crossing the wall of geometric phases of a gauged linear
    sigma model does arise solely from the curve neighborhood of
    $C$ and $C^{\prime}$ that are involved in the topology change
    $X\dashrightarrow X^{\prime}$.
   It is in this way that A.-M.~Li and Y.~Ruan in [L-R]
    completed Witten's picture/wall-crossing-formula in [Wit5]
    mathematically.
 \end{itemize}

\bigskip

\begin{flushleft}
{\large\bf A.2\hspace{1ex}
 Transformation of GW-invariants from the symplectic viewpoint.}
\end{flushleft}
In this appendix, we highlight the work of A.-M.~Li and
  Y.~Ruan [L-R] in the symplectic aspect.
Besides for comparison to the algebraic main text,
 these techniques become important for new applications.
See also [I-P1] and [I-P2] of E.-N.~Ionel and T.H.~Parker.
(Slight modifications of notations are made so that they are
 linked to the main text directly.)

\bigskip

\begin{flushleft}
{\bf Symplectic relative Gromov-Witten theory.}
\end{flushleft}
Let $(M,\partial M; \omega_0)$ be a symplectic manifold-with-boundary
 equipped with a {\it regular local Hamiltonian function}
 $H:U_0\rightarrow (-\infty,0] \subset {\Bbb R}$
 defined on a neighborhood $U_0$ of $\partial M$
  such that $H^{-1}(0)=\partial M$ and
  that the local Hamiltonian flow generated by the Hamiltonian
   vector field $X_H$ on $U_0$ renders $U_0$ a circle bundle
   over the associated quotient manifold $U_0/S^1$ with boundary
   $Z := \partial M /S^1$.
By construction, $Z$ is a symplectic reduction of $U_0$ and
 is equipped with a natural symplectic structure
 $\underline{\omega}$ induced from $(U_0,\omega_0)$.
An {\it infinite symplectic stretching} is then appied to 
 $(M,\omega_0)$ along $\partial M$.
The resulting manifold-with-end is topologically the peeled manifold
 $M^{\circ}:= M-\partial M$ with a collar of $\partial M$ pulled out
 to an infinite end diffeomorphic $[0,+\infty) \times \partial M$.
$M^{\circ}$ is equipped with a symplectic structure $\omega$
 from stretching $\omega_0$ along the end.
One can choose an $\omega$-tamed almost complex structure $J$ on
 $M^{\circ}$ whose restriction to the infinite end is compatible
 with an $\underline{\omega}_0$-tamed almost structure on $Z$ and
 satisfies compatibility conditions with the $S^1$ Hamiltonian
 action on the end.
These conditions imply that when the image of a $J$-holomorphic
 map $f$ enters the end, it has to go all the way through to
 infinity.
Such occurrence corresponds to {\it non-removable singularities}
 of $f$ on the domain Riemann surface.
The pair $(J,\omega)$ defines also a Riemannian metric on $M^{\circ}$
 and on ${\Bbb R}^1\times \partial M$
 that is needed for defining various Sobolev spaces in the problem.

The energy for a $J$-holomorphic map $f$ from a Riemann surface
 $\Sigma$ (with possibly nodes and punctures) to $M^{\circ}$ is
 defined by $E(f):=\int_{\Sigma}f^{\ast}\omega$.
For a {\it $J$-holomorphic map} $f$ to have {\it finite energy},
 any small enough neighborhood of a non-removable singularity of $f$
 on $\Sigma$ must be mapped to the end of $M^{\circ}$ in such a way
 that $f$ is asymptotically wrapping an orbit of the $S^1$-action
 on the end.
Such maps can be identified with a {\it relative map} to
 $(\overline{M^{\circ}},Z)$, where $\overline{M^{\circ}}$
 is the closed manifold obtained from $M$ by collapsing
 the $S^1$-orbits in $\partial M$ via the quotient map
 $\partial M\rightarrow Z$.
The structures $\omega$ and $J$ on $M^{\circ}$ descend naturally
 to the symplectic manifold-submanifold pair
 $(\overline{M^{\circ}},Z)$.

Under deformation of $J$-holomorphic maps, non-vanishing positive
 energy can condense toward isolated points on the domain Riemann
 surface.
This gives rise to the bubbling phenomena of $J$-holomorphic maps.
For a target with infinite ends, there are two types of bubbles that
 can occur:
 (1) the ordinary ones that already occur in the study of absolute
     symplectic Gromov-Witten theory and
 (2) the {\it ghost bubbles}, which are a collection of
     $J$-holomorphic maps of finite enegry from ${\Bbb P}^1$ with
     punctures to ${\Bbb R}^1\times \partial M$.
The latter has to satisfy some compatibility conditions
 for gluing the bubbles along ends at infinity.
(This is the pre-deformability condition in the algebro-geometric
 setting.).
Fix a topological type $(g,T_m;\beta)$ of $J$-holomorphic relative
 maps, where
 $g$ is the arithmetic genus of domain Riemann surface
   (possibly with nodes),
 $T_m$ is an ordered set of $m$ non-negative integers labelling 
   both the ordinary and the special marked points (corresponding
   to non-removable singularities of maps) on the Riemann surface
   and the asymptotic behavior of maps on the special marked points,
   and
 $\beta\in H_2(M^{\circ};{\Bbb Z})$ is a curve class.
Then adding in $J$-holomorphic maps on domains with bubbles
 gives rise to a compact Hausdorff {\it moduli space
 $\overline{\cal M}(M^{\circ}, (g,T_m;\beta))$
 of $J$-holomorphic stable relative maps} of finite energy
 into $M^{\circ}$.
(Such a map can be identified as a stable relative map of
 topological type $(g,\,\cdot\,;\beta)$ to an expanded relative
 pair $(\overline{M^{\circ}}_{[l]}, Z_{[l]})$ from the smooth
 variety-divisor pair $(\overline{M^{\circ}}, Z)\,$.)

In general, $\overline{\cal M}(M^{\circ}, (g,T_m;\beta))$
 may have dimension greater than expected (either from the index
 of a natural differential operator or from the standard
 deformation-obstruction theory in the moduli problem) and
is unsuitable to be directly used to define invariants of
 $(M,\partial M)$ constant under deformations.
The {\it virtual neighborhood construction} in [Ru1],
 outlined below for the current case, can be employed
 to remedy this.
(See [Si1] for a review and references on various alternative
 approaches in the symplectic category.)

Let $\overline{\cal B}(M^{\circ},(g,T_m;\beta))$ be the
 {\it moduli space of $C^{\infty}$ finite-energy stable relative
  maps} of topological type $(g,T_m;\beta)$ without imposing
 the $J$-holomorphy condition: $df+ J\circ df\circ j=0\,$.
Let ${\cal D}_{(g,T_m;\,\beta)}$ be the set of weighted graphs
 with roots and legs that labels the following combinatorial data
 of stable maps:
  (i) the combinatorial/topological type of the domain Riemann
      surface with ordinary and special marked points,
 (ii) the index of the asymptotic wrapping (around the $S^1$-orbits
       in the end of $M^{\circ}$ or ${\Bbb R}^1\times\partial M$)
       associated to each special marked point,
(iii) corresponding decomposition of curve classes
      $\beta=\sum_i\,\beta_i$ associated to the decomposition
      of domain Riemann surface in (i).
Then ${\cal D}_{(g,T_m;\,\beta)}$ is a finite set and
 both $\overline{\cal M}(M^{\circ},(g,T_m;\beta))$ and
 $\overline{\cal B}(M^{\circ},(g,T_m;\beta))$ are naturally
 stratified with each statum labelled by
 a $\Gamma\in {\cal D}_{(g,T_m;\,\beta)}\,$.
Each stratum $\overline{\cal B}(M^{\circ},\Gamma)$ of
 $\overline{\cal B}(M^{\circ},(g,T_m;\beta))$ is a Fr\'{e}chet
 orbifold.
There is an {\it obstruction space fibration}
 $\overline{\cal F}
  \rightarrow \overline{\cal B}(M^{\circ},(g,T_m;\beta))$
 whose restriction to each stratum
  ${\cal F}_{\Gamma}\rightarrow {\cal B}_{\Gamma}$
  is an orbi-bundle of infinite rank.
The fiber of $\overline{\cal F}$ at
 $[\,f:\Sigma\rightarrow
   M^0\, \raisebox{.2ex}{\scriptsize $\coprod$}\:
         ({\Bbb R}^1\times\partial M)\,]$
 comes from $C^{\infty}$-sections of the sheaf 
 $\Omega^{0,1}\otimes f^{\ast}T_{\ast}
    (M^0\, \raisebox{.2ex}{\scriptsize $\coprod$}\:
          ({\Bbb R}^1\times\partial M))$
 on $\Sigma$ with finite appropriate Sobolev norm.
The $J$-holomorphy condition can be realized as a section 
 $\bar{\partial}_J:
   \overline{\cal B}(M^{\circ},(g,T_m;\beta))
   \rightarrow \overline{\cal F}(M^{\circ},(g,T_m;\beta))$
 with
 $\bar{\partial}_J^{-1}(0)
      =\overline{\cal M}(M^{\circ},(g,T_m;\beta))\,$.
Let
 $\pi:{\cal C}\rightarrow \overline{\cal B}(M^{\circ},(g,T_m;\beta))$
 be the universal Riemann surface over
 $\overline{\cal B}(M^{\circ},(g,T_m;\beta))$ and
 $f:{\cal C}\rightarrow
   M^0\, \raisebox{.2ex}{\scriptsize $\coprod$}\:
         ({\Bbb R}^1\times\partial M)$
 be the universal $C^{\infty}$ stable maps.
By passing to a multiple, one can assume that
 the induced symplectic $2$-form on $\overline{M^{\circ}}$,
  still  denoted by $\omega$, is integral.
Let ${\cal L}$ be a complex line bundle on $\overline{M^{\circ}}$
 with a unitary connection such that $c_1({\cal L})=[\omega]$.
Let $\omega_{{\cal C}/{\overline{\cal B}}}$ be the relative
 dualizing sheaf of ${\cal C}$
 over $\overline{\cal B}(M^{\circ},(g,T_m;\beta))$.
Identify finite-energy stable maps to $M^{\circ}$ with relative
 maps to $(\overline{M^{\circ}},\,Z)\,$.
Then, from the push-forward
 $R\pi_{\ast}(\omega_{{\cal C}/\overline{\cal B}}\otimes
             f^{\ast}{\cal L})^{\otimes k}$
 with $k$ large enough, one can construct a {\it stratumwise
 complex orbi-bundle $\rho:\overline{\cal E}\rightarrow {\cal U}$
 of finite constant rank} over a neighborhood ${\cal U}$
 of $\overline{\cal M}(M^{\circ},(g,T_m;\beta))$
 in $\overline{\cal B}(M^{\circ},(g,T_m;\beta))$
such that the section $\bar{\partial}_J$ of $\overline{\cal F}$
 over ${\cal U}$ can be extended to a section ${\cal S}_e$
 of $\rho^{\ast}\overline{\cal F}$ over $\overline{\cal E}$
 with $D{\cal S}_e$ surjective.
This implies the transversality of ${\cal S}_e$ with
 the zero-section of $\rho^{\ast}\overline{\cal F}$ over
 $\overline{\cal E}$.
The data $(\overline{\cal E}, {\cal S}_e)$ is called a
 {\it stabilization} of the section $\bar{\partial}_J$
 of $\overline{\cal F}$ along $\bar{\partial}_J^{-1}(0)$.
Note that $\overline{\cal E}$ over ${\cal U}$ is oriented
 by construction.

The pair $(\overline{\cal E},{\cal S}_e)$ defines 
 a {\it virtual neighborhood} $(U,\,E,\,I)$
 of $\overline{\cal M}(M^{\circ}, (g,T_m;\,\beta))\,$, where
  \begin{itemize}
   \item[$\cdot$]
    $U:= ({\cal S}_e)^{-1}(0)\,$,
     which has dimension the rank of ${\cal E}$ plus the expected
      dimension of the moduli problem;
    $U$ has an induced stratification from that of $\overline{\cal E}$
     labelled by ${\cal D}_{(g,T_m;\,\beta)}\,$,
    this stratification has the property
     that each stratum $U_{\Gamma}$ is a smooth orbifold and
     that
      if ${\cal B}_{\Gamma_1}\subset \overline{{\cal B}_{\Gamma_2}}$
       is a lower stratum,
      then $U_{\Gamma_1}\subset U_{\Gamma_2}$ is a submanifold of
       (real) codimension $\ge 2\,$; 

   \item[$\cdot$]
    the stratumwise orbi-bundle map ${\cal E}\rightarrow {\cal U}$
     induces a map $U\rightarrow {\cal U}$ and
    $E\rightarrow U$ is the pull-back of ${\cal E}$ via this map;
    by construction $E$ is an oriented stratumwise orbi-bundle
     over $U$;

   \item[$\cdot$]
    $I$ is the tautological section of $E$ associated to the inclusion
     $U\subset \overline{\cal E}$;
    by construction
     $I^{-1}(0)=\overline{\cal M}(M^{\circ}, (g,T_m;\,\beta))\,$.
  \end{itemize}
The relative Gromov-Witten invariants can now be defined
 via the virtual neighborhood $(U,E,I)$ of
 $\overline{\cal M}(M^{\circ}, (g,T_m;\beta))$
 as follows.

Let $\Theta$ be the Thom class of the oriented stratumwise
 orbi-bundle $E$ over $U$.
The ordinary and the special marked points on the domain Riemann
 surfaces define eveluation maps from compositions
 $$
  \ev_i\,:\;
    U\; \longrightarrow {\cal M}(M^{\circ}, (g,T_m;\,\beta))\;
        \longrightarrow\;  M^{\circ}
 $$
 and
 $$
  \ev_j\,:\;
    U\; \longrightarrow {\cal M}(M^{\circ}, (g,T_m;\,\beta))\;
        \longrightarrow\;  Z
 $$
 respectively.
Let $\alpha_i\in H^{\ast}(M^{\circ})$ and $\xi_j\in H^{\ast}(Z)$
 represented by differential formas,
then the ({\it symplectic}) {\it relative Gromov-Witten invariants}
 associated to these data is defined to be
 $$
  \Psi^{(M,\,\partial M)}_{(g,T_m;\,\beta)}
   (\alpha_1,\,\cdots\,,\, \xi_1,\,\cdots)\;
  =\; \int_U\, I^{\ast}\Theta\,
               \wedge\, \prod_i\ev_i^{\ast}\,\alpha_i\,
               \wedge\, \prod_j \ev_j^{\ast}\,\xi_j\,.
 $$
The integral is defined stratum by stratum and only the strata of
 maximal dimension matter.
Inequalities that bound appropriate Sobolev norms of
 the normal derivatives of gluing maps of $\ker D{\cal S}_e$
 from one stratum to a lower stratum
imply that the integral on each stratum is finite.
This shows that 
 $\Psi^{(M,\,\partial M)}_{(g,T_m;\,\beta)}
   (\alpha_1,\,\cdots\,,\, \xi_1,\,\cdots)$
 is well-defined.
These will be regarded also as the relative Gromov-Witten invariants
 for the symplectic pair $(\overline{M^{\circ}},Z)$
 and denoted by
 $\Psi^{(\overline{M^{\circ}},\,Z)}_{(g,T_m;\,\beta)}
    (\alpha_1,\,\cdots\,,\, \xi_1,\,\cdots)\,$.


Relative pairs of the form $(\overline{M^0},Z)$ appear in
 symplectic cuts of a symplectic manifold and
the relative Gromov-Witten invariants defined above form
 the basic blocks of a gluing formula of (symplectic) Gromov-Witten
 invariants, which we now turn to.


\bigskip

\begin{flushleft}
{\bf Gluing formula and transformation of GW-invariants.}
\end{flushleft}
Let $M$ be a (closed) symplectic manifold with a regular local
 Hamiltonian function $H$ defined on an open subset $U_0$ of $M$
 such that $H^{-1}(0)$ separates $M$ into a disjoint union
  $M-H^{-1}(0)
   = M^+\, \raisebox{.2ex}{\scriptsize $\coprod$}\, M^-$
 and that the flow of the Hamiltonian vector field $X_H$ generates
  an $S^1$-aaction on $U$, which in particular renders
  $H^{-1}(0)$ a circle bundle $H^{-1}(0)\rightarrow Z$.
Let $\overline{M}^{\pm}$ be the symplectic manifold obtained
 from $M^{\pm}\cup H^{-1}(0)$ by collapsing the $S^1$-orbits
 in $H^{-1}(0)$ via $H^{-1}(0)\rightarrow Z$.
Then $(\overline{M}^+, \overline{M}^-)$ forms a {\it symplectic cut}
 of $M$ and there is a natural morphism
 $\pi: M\rightarrow \overline{M}^+\cup_Z \overline{M}^-$
 that is a symplectomorphism except along $H^{-1}(0)$, where
  $\pi$ sends $H^{-1}(0)$ onto $Z$, becoming the circle bundle map.

An infinite symplectic stretching can be applied to a small
 tubular neighborhood of $H^{-1}(0)$ in the domain of $H$.
The resulting {\it symplectic space with an infinite tube}
 $M_{\infty}$ is a disjoint union of the infinite symplectic
 stretching of $M^+$ and $M^-$ respectively together with
 structure-preserving diffeomorphisms
 $\partial_{\infty}M^+ \simeq \partial_{\infty}M^- \simeq H^{-1}(0)$
 on the ideal boundary.
The same construction for symplectic relative Gromov-Witten
 theory and invariants can be applied to construct
 a {\it Gromov-Witten theory and invariants for $M_{\infty}$}.
{\it Stokes' theorem} then implies that Gromov-Witten invariants
 for $M$ and those for $M_{\infty}$ with matching curve classes
 and cohomology classes are equal.
On the other hand, the Gromov-Witten invariants for $M_{\infty}$
 can be readily expressed as a gluing of relative Gromov-Witten
  invariants of $(\overline{M}^+,Z)$ with relative Gromov-Witten
  invariants of $(\overline{M}^-,Z)$.
Together, this gives a {\it gluing formula} for Gromov-Witten
 invariants of $M$ in terms of relative Gromov-Witten invariants
 of the pairs,
  $(\overline{M}^+,Z)$ and $(\overline{M}^-,Z)$,
([L-R: Theorem 5.7, Theorem~5.8]).

Given {\it $3$-folds} $M$ and $M^{\prime}$ related either by
 a {\it flop} or by a {\it small extremal transition},
one can {\it associate symplectic cuts}
 $\pi:M\rightarrow \overline{M}^+\cup_Z\overline{M}^-$ and
 $\pi^{\prime}:M^{\prime} \rightarrow
    \overline{M^{\prime}}^+\cup_{Z^{\prime}} \overline{M^{\prime}}^-$
 to $M$ and $M^{\prime}$ respectively
 in such a way that
 $(\overline{M}^+,Z)$ and $(\overline{M^{\prime}}^+, Z^{\prime})$
 depends only on the curve neighborhood involved
 in the flop/small-extremal-transition
while $(\overline{M}^-,Z)$ and $(\overline{M^{\prime}}^-,Z^{\prime})$
 are isomorphic.
This enables one to {\it compare Gromov-Witten invariants of both $M$
 and $M^{\prime}$ to relative Gromov-Witten invariants of
 the symplectic pair $(\overline{M}^-,\,Z)$} and leads to
 {\it transformation rules} between Gromov-Witten invariants
 of $M$ and Gromov-Witten invariants of $M^{\prime}$,
[L-R: Theorem A, Corollary A.2, Theorem B, Corollary B.2].

\vspace{4em}

\end{document}

%% file: psfig-97.tex
\catcode`\@=11\relax
\newwrite\@unused
\def\typeout#1{{\let\protect\string\immediate\write\@unused{#1}}}

\def\psglobal#1{
\immediate\special{ps: plotfile #1 }}
\def\psfiginit{\typeout{psfiginit}
\immediate\psglobal{figtex.pro}%
\special{ps:: /TeXMagnification {\the\mag} def}
}

\def\@nnil{\@nil}
\def\@empty{}
\def\@psdonoop#1\@@#2#3{}
\def\@psdo#1:=#2\do#3{\edef\@psdotmp{#2}\ifx\@psdotmp\@empty \else
    \expandafter\@psdoloop#2,\@nil,\@nil\@@#1{#3}\fi}
\def\@psdoloop#1,#2,#3\@@#4#5{\def#4{#1}\ifx #4\@nnil \else
       #5\def#4{#2}\ifx #4\@nnil \else#5\@ipsdoloop #3\@@#4{#5}\fi\fi}
\def\@ipsdoloop#1,#2\@@#3#4{\def#3{#1}\ifx #3\@nnil
       \let\@nextwhile=\@psdonoop \else
      #4\relax\let\@nextwhile=\@ipsdoloop\fi\@nextwhile#2\@@#3{#4}}
\def\@tpsdo#1:=#2\do#3{\xdef\@psdotmp{#2}\ifx\@psdotmp\@empty \else
    \@tpsdoloop#2\@nil\@nil\@@#1{#3}\fi}
\def\@tpsdoloop#1#2\@@#3#4{\def#3{#1}\ifx #3\@nnil
       \let\@nextwhile=\@psdonoop \else
      #4\relax\let\@nextwhile=\@tpsdoloop\fi\@nextwhile#2\@@#3{#4}}
\def\psdraft{
	\def\@psdraft{0}
	\def\@psdraftspecial{100}
}
\def\psdraftspecial{
	\def\@psdraft{0}
	\def\@psdraftspecial{0}
}
\def\psfull{
	\def\@psdraft{100}
}
\psfull

\newif\if@prologfile
\newif\if@postlogfile
\newif\if@bbllx
\newif\if@bblly
\newif\if@bburx
\newif\if@bbury
\newif\if@height
\newif\if@width
\newif\if@rheight
\newif\if@rwidth
\newif\if@clip
\newif\if@right
\newif\if@left
\newif\if@toplines
\newif\if@box
\newif\if@caption
\newif\if@surround
\newif\if@captionwidth
\newif\if@captionwrite
\newif\if@captionopen
\def\@p@@sclip#1{\@cliptrue}
\def\@p@@sfile#1{
		\def\@p@sfile{#1}
}
\def\@p@@sfigure#1{
		\def\@p@sfile{#1}
}
\def\@p@sfake{\hbox to 0pt{\hss Whatever\hss}}
\def\@p@@sbbllx#1{
		\@bbllxtrue
		\@d@mscratch=#1
		\edef\@p@sbbllx{\number\@d@mscratch}
}
\def\@p@@sbblly#1{
		\@bbllytrue
		\@d@mscratch=#1
		\edef\@p@sbblly{\number\@d@mscratch}
}
\def\@p@@sbburx#1{
		\@bburxtrue
		\@d@mscratch=#1
		\edef\@p@sbburx{\number\@d@mscratch}
}
\def\@p@@sbbury#1{
		\@bburytrue
		\@d@mscratch=#1
		\edef\@p@sbbury{\number\@d@mscratch}
}
\def\@p@@sheight#1{
		\@heighttrue
		\@d@mscratch=#1
   		\edef\@p@sheight{\number\@d@mscratch}
}
\def\@p@@swidth#1{
		\@widthtrue
		\@d@mscratch=#1
		\edef\@p@swidth{\number\@d@mscratch}
}
\def\@p@@srheight#1{
		\@rheighttrue
		\@d@mscratch=#1
		\edef\@p@srheight{\number\@d@mscratch}
}
\def\@p@@srwidth#1{
		\@rwidthtrue
		\@d@mscratch=#1
		\edef\@p@srwidth{\number\@d@mscratch}
}
\def\@p@@sright#1{\@righttrue \@surroundtrue}
\def\@p@@sleft#1{\@lefttrue \@surroundtrue}
\def\@p@@sextraheight#1{\@d@mextraheight=#1}
\def\@p@@sbox#1{\@boxtrue}
\def\@p@@scaption#1{\@captiontrue}
\def\@p@@stoplines#1{
		\@toplinestrue
		\@c@ttoplines=#1
}
\def\@p@@scaptionwidth#1{
		\@captionwidthtrue
	  	\@d@mcaptionwidth=#1
}
\def\@p@@scaptionwrite#1{
		\global\@captionwritetrue
		\global\@w@rname=\expandafter{\jobname_captions.tex}
		\typeout{Captions are written to \the\@w@rname}
}

\def\@p@@sprolog#1{\@prologfiletrue\def\@prologfileval{#1}}
\def\@p@@spostlog#1{\@postlogfiletrue\def\@postlogfileval{#1}}
\def\@cs@name#1{\csname #1\endcsname}
\def\@setparms#1=#2,{\@cs@name{@p@@s#1}{#2}}
\def\ps@init@parms{
		\@bbllxfalse \@bbllyfalse
		\@bburxfalse \@bburyfalse
		\@heightfalse \@widthfalse
		\@rheightfalse \@rwidthfalse
		\def\@p@sbbllx{}\def\@p@sbblly{}
		\def\@p@sbburx{}\def\@p@sbbury{}
		\def\@p@sheight{}\def\@p@swidth{}
		\def\@p@srheight{}\def\@p@srwidth{}
		\def\@p@sfile{}
		\def\@p@scost{10}
		\def\@sc{}
		\@prologfilefalse
		\@postlogfilefalse
		\@clipfalse
		\@rightfalse \@leftfalse
		\@boxfalse \@captionfalse
		\@toplinesfalse \@surroundfalse
		\@d@mextraheight=0pt
 		\@c@ttoplines=0
		\@pshape={} \def\@p@srheight@total{}
		\@captionwidthfalse \@d@mcaptionwidth=0pt
}
\def\parse@ps@parms#1{
	 	\@psdo\@psfiga:=#1\do
		   {\expandafter\@setparms\@psfiga,}}
\newif\ifno@bb
\newif\ifnot@eof
\newread\ps@stream
\newtoks\@linetok
\def\bb@missing{
	\typeout{psfig: searching \@p@sfile \space  for bounding box}
	\openin\ps@stream=\@p@sfile
	\no@bbtrue
	\not@eoftrue
	\catcode`\%=12
	\loop
		\read\ps@stream to \line@in
		\global\@linetok=\expandafter{\line@in}
		\ifeof\ps@stream \not@eoffalse \fi
		\@bbtest{\@linetok}
		\if@bbmatch\not@eoffalse\expandafter\bb@cull\the\@linetok\fi
	\ifnot@eof \repeat
	\catcode`\%=14
}	
\catcode`\%=12
\newif\if@bbmatch
\def\@bbtest#1{\expandafter\@a@\the#1
\long\def\@a@#1
     \ifx\@bbtest#2\@bbmatchfalse\else\@bbmatchtrue\fi}
\long\def\bb@cull#1 #2 #3 #4 #5 {
	\@d@mscratch=#2 bp\edef\@p@sbbllx{\number\@d@mscratch}
	\@d@mscratch=#3 bp\edef\@p@sbblly{\number\@d@mscratch}
	\@d@mscratch=#4 bp\edef\@p@sbburx{\number\@d@mscratch}
	\@d@mscratch=#5 bp\edef\@p@sbbury{\number\@d@mscratch}
	\no@bbfalse
}
\catcode`\%=14
\def\compute@bb{
		\no@bbfalse
		\if@bbllx \else \no@bbtrue \fi
		\if@bblly \else \no@bbtrue \fi
		\if@bburx \else \no@bbtrue \fi
		\if@bbury \else \no@bbtrue \fi
		\ifno@bb \bb@missing \fi
		\ifno@bb \typeout{FATAL ERROR: no bb supplied or found}
			\no-bb-error
		\fi
		\count203=\@p@sbburx
		\count204=\@p@sbbury
		\advance\count203 by -\@p@sbbllx
		\advance\count204 by -\@p@sbblly
		\edef\@bbw{\number\count203}
		\edef\@bbh{\number\count204}
}
\def\in@hundreds#1#2#3{\count240=#2 \count241=#3
		     \count100=\count240	
		     \divide\count100 by \count241
		     \count101=\count100
		     \multiply\count101 by \count241
		     \advance\count240 by -\count101
		     \multiply\count240 by 10
		     \count101=\count240	
		     \divide\count101 by \count241
		     \count102=\count101
		     \multiply\count102 by \count241
		     \advance\count240 by -\count102
		     \multiply\count240 by 10
		     \count102=\count240	
		     \divide\count102 by \count241
		     \count200=#1\count205=0
		     \count201=\count200
			\multiply\count201 by \count100
		     	\advance\count205 by \count201
		     \count201=\count200
			\divide\count201 by 10
		     	\multiply\count201 by \count101
			\advance\count205 by \count201
		     \count201=\count200
			\divide\count201 by 100
			\multiply\count201 by \count102
			\advance\count205 by \count201
		     \edef\@result{\number\count205}
}
\def\compute@wfromh{
		\in@hundreds{\@p@sheight}{\@bbw}{\@bbh}
		\edef\@p@swidth{\@result}
}
\def\compute@hfromw{
		\in@hundreds{\@p@swidth}{\@bbh}{\@bbw}
		\edef\@p@sheight{\@result}
}
\def\compute@handw{
		\if@height
			\if@width
			\else
				\compute@wfromh
			\fi
		\else
			\if@width
				\compute@hfromw
			\else
				\edef\@p@sheight{\@bbh}
				\edef\@p@swidth{\@bbw}
			\fi
		\fi
}
\def\compute@resv{
		\if@rheight \else \edef\@p@srheight{\@p@sheight} \fi
		\if@rwidth \else \edef\@p@srwidth{\@p@swidth} \fi
		\edef\@p@srheight@total{\@p@srheight}
}
\newtoks\@pshape
\def\@c@ttoplines{\count120}
\def\@c@theightcount{\count121}
\def\@c@tshapecount{\count122}
\newdimen\@d@mwidthshape
\newdimen\@d@mextraheight
\newdimen\@d@mscratch
\def\compute@parshape{
	\if@right
		\if@left
	   		\typeout{error: Can't have both left and right set}
			\@leftfalse
		\fi
	\fi
	\@d@mscratch=\@p@swidth truesp
	\divide \@d@mscratch by 19
	\multiply \@d@mscratch by 20
	\edef\@p@swidthdimen{\the\@d@mscratch}
	\@c@tshapecount=\@c@ttoplines
 	\@d@mscratch=\baselineskip
	\multiply \@d@mscratch by \@c@ttoplines
	\advance \@d@mscratch by .4\baselineskip
    	\edef\@p@stopdistance{\the\@d@mscratch }
	\@d@mscratch=\@p@sheight truesp
	\divide \@d@mscratch by 2
	\edef\@p@shalfboxheight{\the\@d@mscratch}
	\if@toplines
		\loop \@pshape=\expandafter{\the\@pshape 0pt \hsize}
		\advance\@c@ttoplines by -1
		\ifnum\@c@ttoplines>0 \repeat
	\fi
   	\@c@theightcount=\@p@srheight@total
	\advance \@c@theightcount by \@d@mextraheight
	\divide  \@c@theightcount by \baselineskip
	\advance \@c@theightcount by 1
    	\advance \@c@tshapecount by \@c@theightcount
	\advance \@c@theightcount by -1
	\@d@mwidthshape=\hsize
     	\advance \@d@mwidthshape by -\@p@swidthdimen
	\if@left
		\def\@moveshape{0pt}
		\ifnum\@c@theightcount>0
		      	\loop
			\@pshape=%
			\expandafter{\the\@pshape %
					\@p@swidthdimen \@d@mwidthshape}
			\advance \@c@theightcount by -1
			\ifnum\@c@theightcount>0 \repeat
		\else
			\advance \@c@tshapecount by 1
		\fi
	\fi
	\if@right
		\@d@mscratch=\hsize
		\advance \@d@mscratch by -\@p@swidth truesp
		\edef\@moveshape{\@d@mscratch}
		\ifnum\@c@theightcount>0
			\loop
			\@pshape=\expandafter{\the\@pshape 0pt \@d@mwidthshape}
			\advance \@c@theightcount by -1
			\ifnum\@c@theightcount>0 \repeat
		\else
			\advance \@c@tshapecount by 1
		\fi
	\fi
	\ifnum \@p@srheight=0
		\@pshape={}
		\@c@tshapecount = 0
	\else
	 	\@pshape=\expandafter{\the\@pshape 0pt \hsize}
	\fi
}
\def\@p@ssetsurroundboxes{
	\global\parshape=\count122 \the\@pshape		
 	\moveright\@moveshape
	\vbox to 0pt\bgroup\hskip0pt\vskip\@p@stopdistance
}
\newtoks\@captiontok
\newbox\@b@xcaption
\newdimen\@d@mcaptionwidth
\newdimen\@d@mcaptionheight
\newwrite\@w@rcaption
\newtoks\@w@rname
\def\setcaption#1{\@captiontok={#1}}
\def\@set@caption{
	\setbox\@b@xcaption\vbox{\hsize\@d@mcaptionwidth
	\tolerance=9000 \baselineskip=11.4pt
	\noindent\relax\the\@captiontok}
	\if@captionwrite
		\if@captionopen
		\else
			\global\@captionopentrue
			\immediate\openout\@w@rcaption=\the\@w@rname
		\fi
		\immediate\write\@w@rcaption{\the\@captiontok}
		\immediate\write\@w@rcaption{}
	\fi
}
\def\compute@caption{
	\if@captionwidth
	\else
		\@d@mcaptionwidth = \@p@swidth truesp
		\divide \@d@mcaptionwidth by 20
		\multiply \@d@mcaptionwidth by 17
	\fi
	\@set@caption
	\@d@mcaptionheight=\ht\@b@xcaption
	\if@rheight
	\else
		\count100=\@p@srheight
	   	\advance \count100 by \@d@mcaptionheight
	   	\advance \count100 by \bigskipamount
	   	\advance \count100 by \medskipamount
	   	\edef\@p@srheight@total{\number\count100}
	\fi
}
\newif\if@alreadyjtem \@alreadyjtemfalse
\def\newpar{\hfil\vadjust{\vskip\parskip}%
	{\count100=\parskip
	\count101=\baselineskip
	\divide\count101 by 10  \multiply\count101 by 3
	\advance \count100 by \count101
	\divide\count100 by \baselineskip
	\advance\count100 by \prevgraf
	\global\prevgraf=\count100}%
	\break\if@alreadyjtem\else\indent\fi%
}
\let\sav@par=\par
\def\jtem#1{%
    	\if@alreadyjtem\else\bgroup\fi
	\def\par{\sav@par\egroup\sav@par}
	\if@alreadyjtem\else\leftskip\parindent\fi
	\@alreadyjtemtrue
	\noindent\hskip0pt
	\llap{#1\ }\ignorespaces
}
\def\compute@sizes{%
	\compute@bb
	\compute@handw
  	\compute@resv
	\if@caption
		\compute@caption
	\fi
	\if@surround
		\compute@parshape
	\fi
}
\def\@p@sdospecials{
	\ifnum\@p@scost<\@psdraft
	       	\typeout{psfig: including \@p@sfile \space }
	\fi
	\special{ps::[begin] 	\@p@swidth \space \@p@sheight \space
			\@p@sbbllx \space \@p@sbblly \space
			\@p@sbburx \space \@p@sbbury \space
			startTexFig \space }
	\ifnum\@p@scost<\@psdraft
		\if@clip
			\typeout{(clip)}
			\special{ps:: \@p@sbbllx \space \@p@sbblly \space
				\@p@sbburx \space \@p@sbbury \space
			    	doclip \space }
		\fi
	\fi
	\if@box
		\typeout{(box)}
  		\special{ps:: \@p@sbbllx \space \@p@sbblly \space
			\@p@sbburx \space \@p@sbbury \space
		    	dobox \space }
	\fi
	\ifnum\@p@scost<\@psdraft
		\if@prologfile
	    		\special{ps: plotfile \@prologfileval \space }
		\fi
		\special{ps: plotfile \@p@sfile \space }
    		\if@postlogfile
			\special{ps: plotfile \@postlogfileval \space }
		\fi
	\fi
	\special{ps::[end] endTexFig \space }
}
\newif\if@putinvbox

\def\psfig#1{{%
	\ifhmode%
		\vbox\bgroup
		\@putinvboxtrue
	\else
		\@putinvboxfalse
	\fi
       	\ps@init@parms
	\parse@ps@parms{#1}
       	\compute@sizes
	\if@surround
		\psfig@for@surround
	\else
		\psfig@for@regular
	\fi
	\if@putinvbox
       		\egroup
	\fi
}}
\def\psfig@for@surround{%
	\@p@ssetsurroundboxes
	\ifnum\@p@scost<\@psdraft
		\@p@sdospecials
		\vbox to \@p@srheight true sp{\vss}
       	\else
		\if@box
			\@p@sdospecials
		\fi
		\vbox to \@p@srheight true sp{
			\vskip\@p@shalfboxheight
			\hbox to \@p@srwidth true sp{
				\hss
				\ifnum\@p@scost<\@psdraftspecial
					\@p@sfile
				\else
					\@p@sfake
				\fi
      				\hss
			}
		\vss
		}
	\fi
	\if@caption
		\medskip
		\hbox to \@p@srwidth true sp{
			\hss
			\box\@b@xcaption
			\hss
		}
 		\medskip
	\fi
	\vss\egroup
	\vskip-\parskip
}

\def\psfig@for@regular{%
	\if@putinvbox
	\else
		\vskip\parskip
	\fi
	\ifnum\@p@scost<\@psdraft
		\@p@sdospecials
		\vbox to \@p@srheight true sp{%
			\hbox to \@p@srwidth true sp{
			\hfil
			}
		\vfil
		}
       	\else
		\if@box
			\@p@sdospecials
		\fi
	    	\vbox to \@p@srheight true sp{
			\vss
			\hbox to \@p@srwidth true sp{
				\hss
				\ifnum\@p@scost<\@psdraftspecial
					\@p@sfile
				\else
					\@p@sfake
				\fi
				\hss
			}
		    	\vss
		}
	\fi
	\if@caption
		\medskip
		\hbox to \@p@srwidth true sp{
			\hss
			\box\@b@xcaption
			\hss
		}
		\bigskip
	\fi
	\if@putinvbox
	\else
		\vskip-\parskip
	\fi
}
\catcode`\@=12\relax

%% file: newcmd.tex
\newcommand{\Amodel}{\mbox{\it A-model}\,}

\newcommand{\Aut}{\mbox{\it Aut}\,}
\newcommand{\Bl}{\mbox{\it Bl}}

\newcommand{\Eq}{\mbox{\rm Eq}}

\newcommand{\GW}{\mbox{\it GW}\,}

\newcommand{\Spec}{\mbox{\it Spec}\,}

\newcommand{\Spin}{\mbox{\it Spin}\,}

\newcommand{\String}{\mbox{\it String}\,}
\newcommand{\SU}{\mbox{\it SU}\,}

\newcommand{\Wilson}{\mbox{\it\tiny Wilson}}

\newcommand{\IIA}{\mbox{\it IIA}\,}
\newcommand{\IIB}{\mbox{\it IIB}\,}

\newcommand{\bboxtimes}{\Box\hspace{-1.75ex}\raisebox{.15ex}{$\times$}\,}

\newcommand{\dimm}{\mbox{\it dim}\,}
 \newcommand{\codimm}{\mbox{\it codim}\,}
\newcommand{\ev}{\mbox{\it ev}\,}

\newcommand{\longrightaarrow}{\longrightarrow\hspace{-3ex}\longrightarrow}

\newcommand{\rel}{\mbox{\it\scriptsize rel}}

\newcommand{\sing}{\mbox{\scriptsize \it sing}}

\newcommand{\vdim}{\mbox{\it vdim}\,}
\newcommand{\virt}{\mbox{\scriptsize\it virt}}
